\documentclass[preprint,12pt]{elsarticle}

\usepackage{amsfonts}
\usepackage{cite,url,subfigure,epsfig,graphicx}
\usepackage{amssymb,amsmath}
\usepackage{enumitem}
\usepackage{calc}
\usepackage[margin=1in]{geometry}
\usepackage[T1]{fontenc}

\journal{Applied Energy}

\begin{document}

\begin{frontmatter}

\title{Optimal Energy Management for Commercial Buildings Considering Comprehensive Comfort Levels in a Retail Electricity Market}

\author[add1,add2]{Zheming~Liang}
\author[add1]{Desong~Bian}
\author[add1]{Xiaohu~Zhang}
\author[add1]{Di~Shi\corref{cor1}}
\ead{di.shi@geirina.net}
\author[add1]{Ruisheng~Diao}
\author[add1]{Zhiwei~Wang}

\cortext[cor1]{Di Shi is the corresponding author.}
\address[add1]{GEIRI North America, San Jose, CA, USA}
\address[add2]{Department of Electrical and Computer Engineering, University of Michigan-Dearborn, MI, USA}
\fntext[label]{This work is funded by SGCC Hybrid Energy Storage Management Platform for Integrated Energy System.}

\begin{abstract}
Demand response has been implemented by distribution system operators to reduce peak demand and mitigate contingency issues on distribution lines and substations. Specifically, the campus-based commercial buildings make the major contributions to peak demand in a distribution system. Note that prior works neglect the consumers' comfort level in performing demand response, which limits their applications as the incentives are not worth as compared to the loss in comfort levels for most time. Thus, a framework to comprehensively consider both operating costs and comfort levels is necessary. Moreover, distributed energy resources are widely deployed in commercial buildings such as roof-top solar panels, plug-in electric vehicles, and energy storage units, which bring various uncertainties to the distribution systems, i.e., (i) output of renewables; (ii) electricity prices; (iii) arrival and departure of plug-in electric vehicles; (iv) business hour demand response signals and (v) flexible energy demand. In this paper, we propose an optimal demand response framework to enable local control of demand-side appliances that are usually too small to participate in a retail electricity market. Several typical small demand side appliances, i.e., heating, ventilation, and air conditioning systems, electric water heaters and plug-in electric vehicles, are considered in our proposed model. Their operations are coordinated by a central controller, whose objective is to minimize the total cost and maximize the customers' comfort levels for multiple commercial buildings. A scenario-based stochastic programming is leveraged to handle the aforementioned uncertainties. Numerical results based on the practical data demonstrate the effectiveness of the proposed framework. In addition, the trade-off between the operation costs of commercial buildings and customers' comfort levels is illustrated.
\end{abstract}


\end{frontmatter}


\section*{Nomenclature}

\subsection*{Indices and Sets}
\begin{description}[leftmargin=!,labelwidth=\widthof{\bfseries $g_{\text{RT}, t}^{s, +}$, $g_{\text{RT}, t}^{s, -}$}]
\item[$i$] Index of CB.
\item[$k$] Index of bus node.
\item[$t$] Index of time slot.
\item[$j$] Index of EWH.
\item[$n$] Index of ES unit.
\item[$s$] Index of scenario.
\item[$m$] Index of renewable.
\item[$N_i$] Number of CBs.
\item[$N_t$] Number of time slots.
\item[$N_j$] Number of EWHs.
\item[$N_k$] Number of PEVs.
\item[$N_n$] Number of ES units.
\item[$N_s$] Number of scenarios.
\item[$N_m$] Number of renewables.
\end{description}

\subsection*{Parameters}
\begin{description}[leftmargin=!,labelwidth=\widthof{\bfseries $g_{\text{RT}, t}^{s, +}$, $g_{\text{RT}, t}^{s, -}$}]
\item[$c_k$] Degradation cost of $k$-th PEV ($\$ / kWh$).
\item[$c_n$] Degradation cost of $n$-th ES unit ($\$ / kWh$).
\item[$c_p$] Penalty cost ($\$/kW$).
\item[$\zeta_j$] Power to heat ratio.
\item[$C_{\text{water}}$] Specific heat capacity ($J/(kg ^\circ C)$).
\item[$c_t^b$] Gas price at time $t$ ($\$ / kBtu$).
\item[$\overline{H}_j$] Upper bound of heat generation for $j$-th heat boiler $kBtu$.
\item[$L_j^s$] Total energy usage allowed for EWH $j$ ($kJ$).
\item[$M_j$] Mass of water in tank $j$ ($kg$).
\item[$D_t^s$] Aggregated real power demand ($kW$).
\item[$Q_t^s$] Aggregated heat demand ($kBtu$).
\item[$T_j^a, T_j^d$] Begin time and end time of the EWH $j$ ($h$).
\item[$c_{\text{DA}, t}^{s, +}$] Day-ahead electricity buying price at time $t$ ($\$ / kWh$).
\item[$c_{\text{DA}, t}^{s, -}$] Day-ahead electricity selling price at time $t$ ($\$ / kWh$).
\item[$c_{\text{RT}, t}^{s} $] Real-time electricity price at time $t$ ($\$ / kWh$).
\item[$\rho_s$] Probability of scenario $s$.
\item[$w_{m, t}^s$] Output of renewable $m$ at time $t$ ($kW$).
\item[$E_{k}^d$] Desired charging status of PEV $k$ ($kWh$).
\item[$\eta_{k}^{+}$] Charging efficiency of PEV $k$.
\item[$P_{k}^{+}$] Upper bound of charging power for PEV $k$ ($kW$).
\item[$\underline{E}_k$,$\overline{E}_k$] Lower bound and upper bound of PEV $k$ ($kWh$).
\item[$E_{k}^{\text{base}}$] Base energy required for the $k$-th PEV with round trip between its home and the CB system ($kWh$).
\item[$I_{k, t}^s$] Binary PEV arrival and departure indicator at time $t$.
\item[$C_{j}^{\text{min}}$] Lower bound of the water temperature ($^\circ C$).
\item[$C_j^d$] Desired water temperature ($^\circ C$).
\item[$\delta_{j}$] Maximum allowed water temperature deviation ($^\circ C$).
\item[$\underline{l}_j,\overline{l}_j$] Lower bound and upper bound of electricity usage of EWH $j$ ($kW$).
\item[$\eta_{n}^{+},\eta_n^-$] Charging and discharging efficiency of ES unit $n$.
\item[$\underline{E}_n$,$\overline{E}_n$] Lower bound and upper bound of energy level for ES unit $n$ ($kWh$).
\item[$P_{n}^{+}$,$P_n^-$] Upper bound of charging/discharging power for ES unit $n$ ($kW$).
\item[$C_{i}^{\text{max}},C_{i}^{\text{min}}$] Upper bound and lower bound of the indoor temperature ($^\circ C$).
\item[$C_t^{s, \text{out}}$] Outdoor temperature at time $t$ ($^\circ C$).
\item[$\Psi_t^s$] Solar irradiance at time $t$ ($kW/m^2$).
\item[$\sigma_{i, t}$] Cooling/Heating indicator at time $t$.
\item[$\alpha_i$,$\beta_i$] Temperature coefficients.
\item[$\eta_i$] Performance coefficient of HVAC system.
\item[$C_i^d$] Desired indoor temperature ($^\circ C$).
\item[$\delta_i$] Maximum allowed indoor temperature deviation ($^\circ C$).
\end{description}

\subsection*{Variables}
\begin{description}[leftmargin=!,labelwidth=\widthof{\bfseries $g_{\text{RT}, t}^{s, +}$, $g_{\text{RT}, t}^{s, -}$}]
\item[$C_{i,t}^{s, \text{in}}$] Indoor temperature of the $i$-th CB at time $t$ ($^\circ C$).
\item[$C_{i,t}^{s, \text{iw}}$] Inner wall temperature of the $i$-th CB at time $t$ ($^\circ C$).
\item[$C_{i,t}^{s, \text{ow}}$] Outer wall temperature of the $i$-th CB at time $t$ ($^\circ C$).
\item[$P_{i, t}^{s, \text{hvac}, o}$] Output power of the HVAC system at time $t$ ($kW$).
\item[$P_{i, t}^{s, \text{hvac}}$]  Power supplied to the HVAC system at time $t$ ($kW$).
\item[$J_{i,t}^s$] Consumer's comfort level related with the $i$-th HVAC system at time $t$.
\item[$E_{k, t}^s$] Energy level of PEV at time $t$ ($kWh$).
\item[$p_{k, t}^{s, +}$] Power charged into the $k$-th PEV at time $t$ ($kW$).
\item[$J_{k,t}^s$] Comfort level related with the $k$-th PEV at time $t$.
\item[$E_{n, t}^s$] Energy level of ES unit at time $t$ ($kWh$).
\item[$p_{n, t}^{s, +}$, $p_{n, t}^{s, -}$] Power charged into or discharged from the $n$-th ES unit at time $t$($kW$).
\item[$u_{n, t}^{s, +}$, $u_{n, t}^{s, -}$] Binary variables indicating the charging and discharging decisions of the $n$-th ES unit at time $t$.
\item[$h_{j, t}^s$] Heat generation of the $j$-th heat boiler at time $t$ ($kBtu$).
\item[$g_{\text{RT}, t}^{s, +}$, $g_{\text{RT}, t}^{s, -}$] Real-time power buy from and sell to the retail electricity market at time $t$ ($kW$).
\item[$g_{\text{DA}, t}^{-}$, $g_{\text{DA}, t}^{+}$] Day-ahead power sell to and buy from the retail electricity market at time $t$ ($kW$).
\item[$\psi_t$] Auxiliary variable.
\item[$l_{j, t}^s$] Power delivered to the $j$-th EWH at time $t$ ($kW$).
\item[$C_{j, \tau}^s$] Water temperature in the $j$-th EWH ($^\circ C$).
\item[$J_{j,\tau}^s$] Comfort level related with $j$-th EWH.
\end{description}

\section{Introduction}\label{sec:introduction}
Driven by the goal of reducing greenhouse gas emissions and fuel costs, plug-in electric vehicles (PEVs) have attracted lots of attention in recent years. Nevertheless, recent research has shown that even a modest penetration of PEVs could cause frequent overloads in a distribution system if they are allowed to charge whenever they are plugged into the system~\citep{SaOr15}. To avoid such problems, demand response (DR) is introduced as an effective approach to reducing peak demand. DR can be minimum operating cost-driven~\citep{LiAl17}, least load curtailment-driven~\citep{BiPi15}, or comfort-driven~\citep{KoBa16}. However, the combination of these aspects has not been fully investigated. In addition, most of the previous works focus on the DR problems for residential households~\citep{Ra16} rather than commercial buildings (CBs).

Power consumption in buildings accounts for approximately $40\%$ of the total energy usage all over the world~\citep{KlKw12}. In the U.S., about half of the buildings' energy usage are consumed by CBs~\citep{CB}. In a distribution system, a joint operation of CBs would have more advantages than independent operation of individual CB for the following reasons: (i) the load demand of multiple CBs in a campus is larger than that of individual CB, thus the campus-based CBs can act as a large electricity consumer; and (ii) surplus power of CBs can be sold back to the grid, or be shared with other CBs. Generally, large consumers have the privilege to directly negotiate with utilities, and are more likely to obtain favorable electricity prices. Each CB can be treated as a single controllable entity, which can easily participate in a retail electricity market and perform DR functions~\citep{Ra162}. Given these considerations, an efficient energy management framework, which directly works for the commercial campus, should be developed to optimally coordinate the energy scheduling for multiple CBs.

Some prior works neglect consumers' comfort in performing DR~\citep{Ra15}, which limits their applications as the incentives are not worth as compared to the loss in comfort levels for most cases. In addition, only the comfort level of heating, ventilation, and air conditioning (HVAC) system is considered for the previous comfort-driven DR works~\citep{NgLe14Comfort}. The comfort levels regarding PEVs and electric water heaters (EWHs) are not addressed. Thus, a model to simultaneously optimize the operating costs and comprehensive comfort levels is necessary. Moreover, to develop an efficient energy scheduling model incorporating DR process, several uncertainties should be considered, ie., (i) output of renewables; (ii) electricity prices; (iii) arrival and departure of PEVs; (iv) DR signals during the business hours; and (v) flexible energy demand. A customized optimization algorithm should be designed to handle these uncertainties together.

Several optimization approaches have been implemented in prior works to handle only parts of the aforementioned uncertainties. A deterministic approach is presented in~\citep{HuRh17}, while demand flexibility and comfort are quantified for cooling and heating of CBs in hot/cold climate zone. They find out that under a warmer climate, DR appliances offer higher but shorter flexibility, and vice versa. Even though the deterministic approach converges fast, its optimality highly depends on an accurate forecast of all uncertainties~\citep{ChXu17}, which is difficult to obtain. The genetic algorithm (GA) is introduced in~\citep{CuGa17} to save both cooling and heating storages' life-cycle cost through the DR process in CBs. However, GA is limited by its time-consuming convergence iteration and the number of population restricts its optimality, which is not suitable for a fast DR approach. Scenario-based stochastic programming is employed in~\citep{KiNo17} to provide the optimal price-driven DR for CBs with thermal storage units. The uncertainties related to outputs of renewables, day-ahead and real-time electricity prices are considered. A multi-objective DR scheduling for a system with high penetration of renewables is studied in~\citep{HaSh18}; however, comfort levels are not considered. The authors in~\citep{DoMo17} propose a general stochastic optimization framework for the energy scheduling in CBs considering wind power uncertainties. However, the uncertainties from DR appliances, electricity prices and PEVs are not addressed. In summary, none of the aforementioned works has fully considered all distinct factors of DR in commercial campuses.

In this paper, we propose a novel two-stage scenario-based stochastic programming framework to handle uncertainties from DR process in a commercial campus. Our objective is to perform optimal DR to simultaneously minimize the total operating cost of CBs and maximize customers' comfort levels. Comfort levels are not integrated with cost functions using weighting factors; instead, operating costs are unified to follow the same scales as comfort levels. A Monte Carlo method is introduced to generate scenarios for all uncertainties based on historical data. A scenario reduction method is leveraged to reduce a huge number of scenarios to a reasonable number to make the problem tractable. Extensive simulation results based on real-world data sets indicate that our proposed framework is able to minimize operating costs of a commercial campus while maximizing consumers' comfort levels.

In summary, the main contributions of this paper are as follows.
\begin{enumerate}
    \item Several novel piece-wise linear models for comfort levels are proposed to capture unique features of major components in distributed CBs, such as HVAC system, PEVs, and EWHs.

    \item A stochastic scenario-based co-optimization method is introduced to minimize expected operating costs (considering both gas \& electricity) and maximize comprehensive comfort levels simultaneously.

    \item Comfort levels are not transformed into costs with simple weighting factors; instead, we unify the operating costs to follow the same scales with comfort levels.

    \item Extensive simulation results show the modeling of comfort levels is suitable in various scenarios, while our DR framework is reliable to handle all related uncertainties considering comfort levels.
\end{enumerate}

The rest of this paper is organized as follows. In Section~\ref{sec:model}, the models of commercial campus are described. Section~\ref{sec:problem} presents the mathematical formulation of our proposed two-stage scenario based stochastic optimization problem. In Section~\ref{sec:solution}, an efficient solution approach to reformulate and solve the optimization problem is presented. We conduct a case study in Section~\ref{sec:simulation} and then draw conclusions in Section~\ref{sec:conclusion}.

\section{System Modeling}\label{sec:model}

\subsection{Commercial Campus}
\begin{figure}[!ht]
\centering
\includegraphics[width=4.8in]{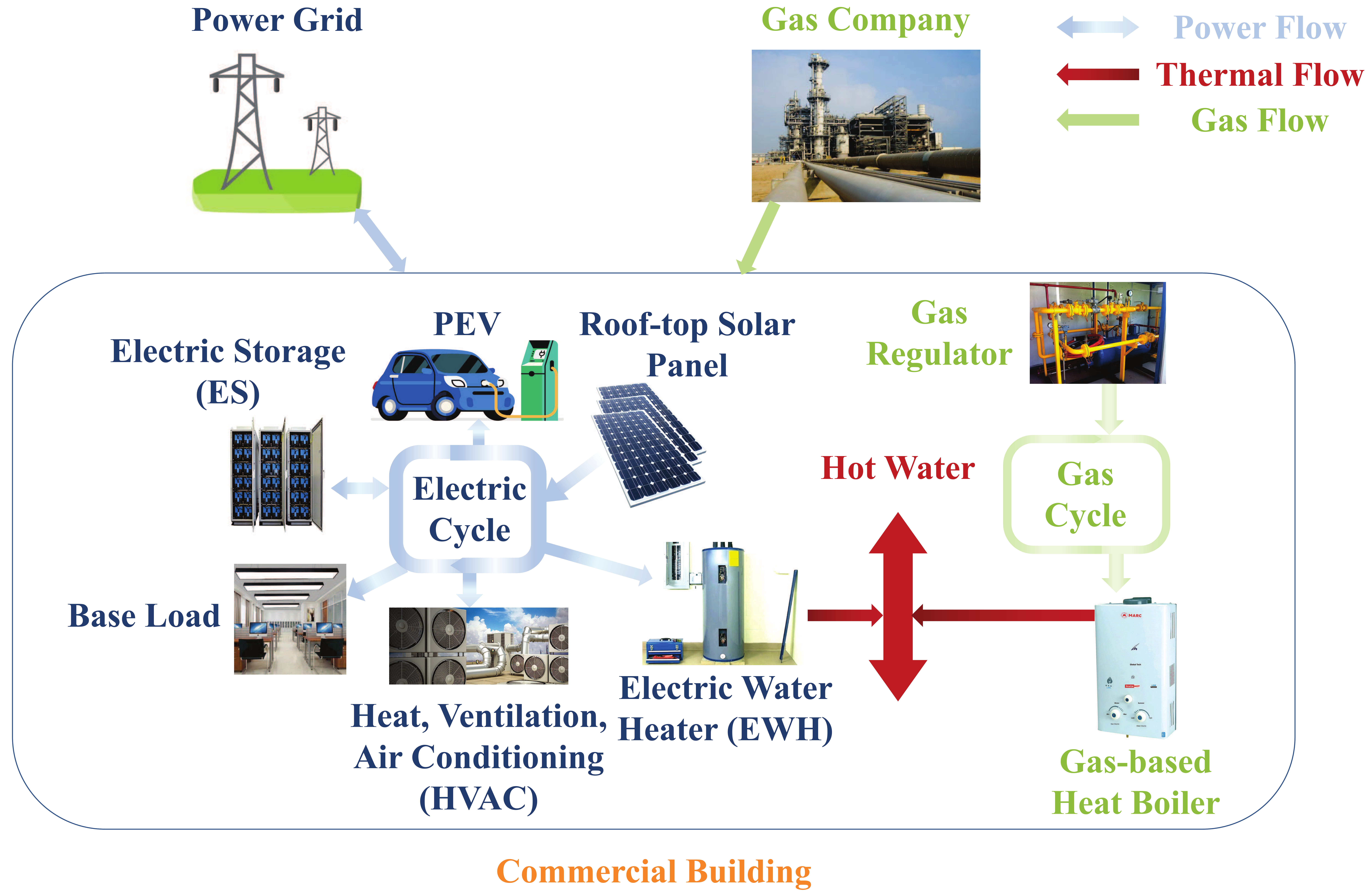}
\caption{A commercial campus with major components.}
\label{fig:System}
\end{figure}
In our model, a commercial campus consists of several CBs, several parking lots with PEV charging stations, and renewables, such as roof-top solar panels. Specifically, for each CB, it includes electricity-based appliances and gas-based appliances, as shown in Fig.~\ref{fig:System}.

For electricity-based appliances, there are several HVAC systems, several electrical storage (ES) units, several EWHs, and several power loads. HVAC systems are installed to provide sufficient cooling airs for consumers inside a CB (in summer). ES units, such as battery packs, are utilized to mitigate power imbalances caused by uncertain output of renewables. EWHs are used to provide hot water for occupants in a CB. Power loads are classified into two categories, i.e., critical power loads and deferrable power loads. Critical power loads (also known as base loads), such as servers, lights and personal computers in a CB, must be satisfied in a whole operating day. Each operating day includes $24$ hours, while each hour is separated into $4$ time slots. Deferrable power loads are both interruptible and shiftable, such as EWHs and PEVs, which have a total amount of energy to be delivered within certain time intervals.

For gas-based appliances, there are several heat boilers and several heat loads. Gas-based heat boilers provide most of the hot water for consumers in a CB when gas prices are much lower than electricity prices or EWHs are not enough. Heat loads are critical loads for hot water, which is provided by both heat boilers and EWHs. This is where electric circle and gas circle on a commercial campus are connected together.

Each appliance in our system has a local controller (LC) that can help the central controller of CBs to perform fast demand response. The commercial campus is connected downstream of both a distribution system operated by a DSO and a gas regulator operated by a gas company. Additionally, the DSO also operates a two-settlement pool-based retail electricity market that has a day-ahead financial market and a real-time physical market~\citep{LiSu18}. The central controller of CBs bids in the day-ahead electricity market based on the forecast of all uncertainties, before an operating day. Then the DSO clears the day-ahead electricity market with day-ahead electricity prices. During the operating day, the DSO clears the real-time electricity market every $15$ minutes based on the real-time power exchanges between commercial campus and distribution system. A penalty would occur when the real-time power exchanges and day-ahead bidding schedules are not the same. Gas prices are based on a long-term contract which is more stable compared with electricity prices.

\subsubsection{HVAC System}
A HVAC system includes heating, ventilation, and cooling functions which can be regarded as a reliable DR appliance~\citep{Ra16}. Additionally, since under most situations, a HVAC system in CBs provide either cooling (summer) or heating (winter). Therefore, we model the ON/OFF status of a HVAC system based on indoor temperatures~\citep{AjLu17}. We model HVAC's cooling/heating dynamics with the third order state-space equations as follows:
\begin{equation}\label{con:HVAC1}
C^s_{i, t+1} = \beta_i C^s_{i, t} + \alpha_i V^s_{i, t}, C^{s, \text{in}}_{i, t} = \Gamma C^s_{i, t}, \forall i, t, s,
\end{equation}
where $C_{i,t}^s = [C_{i,t}^{s, \text{in}}, C_{i,t}^{s, \text{iw}}, C_{i,t}^{s, \text{ow}}]^T$ denotes the state vector. Specifically, $C_{i,t}^{s, \text{in}}$ is the indoor temperature of the $i$-th CB ($^\circ C$). The $C_{i,t}^{s, \text{iw}}$ represents the inner wall temperature of the $i$-th CB ($^\circ C$). The $C_{i,t}^{s, \text{ow}}$ denotes the outer wall temperature of the $i$-th CB ($^\circ C$). $V_{i ,t}^s = [C_{t}^{s, \text{out}}, \Psi_t^s, \sigma_{i, t} P_{i, t}^{s, \text{hvac}, o}]^T$ indicates the input control vector, where $P_{i, t}^{s, \text{hvac}, o} = \eta_i P_{i, t}^{s, \text{hvac}}$. $C_t^{s, \text{out}}$ is the outdoor temperature at time $t$, $\Psi_t^s$ represents the solar irradiance at time $t$, $\sigma_{i, t}$ is the cooling/heating indicator with ($\pm 1$), and $P_{i, t}^{s, \text{hvac}, o}$ is the output power of a HVAC system. $\eta_i$ denotes the coefficient of HVAC system's performance and $P_{i, t}^{s, \text{hvac}}$ represents the power supplied to a HVAC system. $\alpha_i$ and $\beta_i$ are coefficients of CB $i$ that can be derived from the effective window area, the fraction of solar irradiation entering inner walls, the thermal capacitance, and the thermal resistance data. Finally, $\Gamma$ is equal to [1,0,0], which relates the indoor temperature to the state vector~\citep{ZoKu12}.

The indoor temperature of a building is constrained by~\eqref{con:HVAC2} and the output power of a HVAC is bounded in constraint~\eqref{con:HVAC3}.
\begin{equation}\label{con:HVAC2}
C^d_{i} - \delta_{i} \leq C_{i, t}^{s, \text{in}} \leq C^d_{i} + \delta_{i}, \forall i, t, s,
\end{equation}
where $C^d_{i}$ is the desired indoor temperature. $\delta_{i}$ is the maximum allowed temperature deviation from the desired indoor temperature.
\begin{equation}\label{con:HVAC3}
0 \leq P_{i, t}^{s, \text{hvac}} \leq \overline{P}_{i}^{\text{hvac}}, \forall i, t, s,
\end{equation}
where $\overline{P}_{i}^{\text{hvac}}$ is the maximum power consumption of a HVAC system.

The consumers' comfort level related to a HVAC system in the $i$-th CB can be defined as:
\begin{equation}
J_{i,t}^s = \begin{cases}\label{con:comfort1}
0, & C_{i, t}^{s, \text{in}} \geq C_{i}^{\text{max}}, \\
\frac{C_{i}^{\text{max}} - C_{i, t}^{s, \text{in}}}{\delta_{i} - \epsilon_{i}}, & C^d_{i} + \epsilon_{i} \leq C_{i, t}^{s, \text{in}} \leq C_{i}^{\text{max}}, \\
1, & C^d_{i} - \epsilon_{i} \leq C_{i, t}^{s, \text{in}} \leq C^d_{i} + \epsilon_{i}, \\
\frac{C_{i, t}^{s, \text{in}} - C_{i}^{\text{min}}}{\delta_{i} - \epsilon_{i}}, & C_{i}^{\text{min}} \leq C_{i, t}^{s, \text{in}} \leq C^d_{i} - \epsilon_{i},\\
0, & C_{i, t}^{s, \text{in}} \leq C_{i}^{\text{min}}.
\end{cases}
\end{equation}
In equation~\eqref{con:comfort1}, $0$ is the most uncomfort situation, while $1$ is the most comfort one. The comfort indoor temperature zone can be defined as $C^d_{i} \pm \epsilon_{i}$, where $\epsilon_{i}$ is the maximum indoor temperature deviation from the desired temperature that can still ensure a comfort temperature zone. $C_{i}^{\text{max}} = C^d_{i} + \delta_{i}$ and $C_{i}^{\text{min}} = C^d_{i} - \delta_{i}$ are the maximum and minimum tolerable indoor temperatures, respectively.

\subsubsection{PEVs}
PEVs are flexible DR components that have been vastly deployed in many states, especially in California. Besides, charging piles are installed in CBs' parking lots, offering free charging (limited hours) to attract and benefit PEV drivers. However, as mentioned in Section~\ref{sec:introduction}, even a modest penetration of PEVs could cause frequent overloads in a distribution system if they are allowed to charge whenever they are plugged into a system. Therefore, it is important to model their intermittent charging process, while arrival and departure time of PEVs are uncertain~\citep{ZeVa18}.

We use $E_{k, t}^s$ to denote energy stored in the $k$-th PEV at the end of period $t$, with initial available energy $E_{k, 0}^s$. Then, we have the following charging dynamics constraint:
\begin{equation}\label{con:pev_dynamics1}
E_{k, t}^s = E_{k, t - 1}^s + p_{k, t}^{s, +} \eta_{k}^{+} t, \forall k, t, s,
\end{equation}
where $p_{k, t}^{s, +}$ is the power charged into the $k$-th PEV at time $t$ in scenario $s$, and $\eta_{k}^{+}$ represents the charging efficiency of the $k$-th PEV.

Each PEV has a finite capacity, therefore, energy stored in it must have the following lower and upper bounds:
\begin{equation}\label{con:pev_dynamics2}
\underline{E}_k \leq E_{k, t}^s \leq \overline{E}_{k}, \forall k, t, s,
\end{equation}
where the upper bound $\overline{E}_{k}$ is the storage capacity for the $k$-th PEV and the lower bound $\underline{E}_k$ is imposed to reduce the impacts of deep discharging on storage lifetimes (e.g., 5\% of its capacity for a Li-ion battery).

Moreover, the charging rate limits of PEVs are described as follows:
\begin{equation}\label{con:pev_dynamics3}
0 \leq p_{k, t}^{s, +} \leq P_{k}^{+} I_{k,t}^s, \forall k, t, s,
\end{equation}
where $P_k^{+}$ denotes the maximum charged energy over period $t$ for the $k$-th PEV. We use an auxiliary variable $I_{k,t}^s$ to represent the uncertain parking conditions. $I_{k,t}^s=1$ denotes that the $k$-th PEV is parked at a charging station and is available to be charged. On the other hand, $I_{k,t}^s=0$ represents that the $k$-th PEV is not able to be charged.

Furthermore, the comfort level related to the $k$-th PEV can be defined as follows:
\begin{equation}
J_{k,t}^s = \begin{cases}
1, & E_{k}^d \leq E_{k,t}^s, \\ \label{con:comfort2}
\frac{E_{k, t}^s - E_{k}^{\text{base}}}{E^d_{k} - E_{k}^{\text{base}}}, & E_{k}^{\text{base}} \leq E_{k, t}^s \leq E_{k}^d, \\
0, & E_{k, t}^{s} \leq E_{k}^{\text{base}}.
\end{cases}
\end{equation}
In equation~\eqref{con:comfort2}, $0$ is the most uncomfort situation, while $1$ is the most comfort one. The $J_{k,t}^s$ denotes the comfort level of the $k$-th PEV owner. The $E_{k}^d$ is the desired energy stored in the $k$-th PEV. The $E_{k}^{\text{base}}$ represents the base energy required for a round trip of the $k$-th PEV between its home and a CB.

\subsubsection{Electric Water Heater}
For EWHs, they only require a certain amount of electric energy over a specified time interval to serve part of the total heating demand and therefore, have some flexibility in their power profiles. The total electric energy $L_j^s$ for the $j$-th EWH is required to be served within a time interval $[T_j^{a}, T_j^{d}]$, with a minimum and a maximum load serving rate of $\underline{l}_j$ and $\overline{l}_j$, respectively. The $\zeta_j$ denotes the power-to-heat ratio of EWH $j$. The constraints for EWHs are expressed as:
\begin{align}\label{con:EWH1}
& \sum_{t = T_j^{a}}^{T_j^{d}} l_{j, t}^s = L_j^s, \underline{l}_j \leq l_{j, t}^s \leq \overline{l}_j, \forall j, t \in [T_j^a, T_j^d], s \\
& l_{j, t}^s = 0, \forall j, t \notin [T_j^a, T_j^d], s,
\end{align}
where $l_{j, t}^s$ is the power delivered to the $j$-th EWH over period $t$ for scenario $s$.

Moreover, we set the initial water temperature in the $j$-th EWH as $C_{j,0}^{s}$. Thus, the water temperature in the $j$-th EWH can be defined as:
\begin{equation}\label{con:EWH3}
C_{j, \tau}^s = C_{j,0}^s + \Delta C_j^s, \Delta C_j^s = \sum_{t=1}^{\tau} \frac{\zeta_j l_{j, t}^s + h_{j, t}^s - H_{j, t}^{\text{de}}}{M_j C_{\text{water}}}, \forall j, \tau, s,
\end{equation}
where $\Delta C_j^s$ is the temperature deviation of the $j$-th EWH in each time slot, $H_{j, t}^{\text{de}}$ is the heat decrease of the $j$-th EWH in time slot $t$, including the heat loss transferred to its ambient, outflow of the hot water and inflow of the cold water. $M_j$ is the mass of water in tank $j$, and $C_{\text{water}}$ is the specific heat capacity of water. As hot water is supported by both EWHs and gas-based heat boilers, $h_{j, t}^s$ is the hot water provided by the $j$-th heat boiler at time $t$.

Furthermore, the water temperature in the $j$-th EWH needs to be regulated within a certain range to maintain the comfort level related with EWHs:
\begin{equation}\label{con:EWH4}
C_{j}^d - \delta_{j} \leq C_{j, \tau}^s \leq C_{j}^d + \delta_{j}, \forall j, \tau, s
\end{equation}
\begin{equation}
J_{j,\tau}^s = \begin{cases}\label{con:comfort3}
1, & C_{j}^d \leq C_{j, \tau}^s, \\
\frac{C_{j, \tau}^s - C_{j}^{\text{min}}}{C^d_{j} - C_{j}^{\text{min}}}, & C_{j}^{\text{min}} \leq C_{j, \tau}^s \leq C_{j}^d, \\
0, & C_{j, \tau}^{s} \leq C_{j}^{\text{min}}.
\end{cases}
\end{equation}
In equation~\eqref{con:comfort3}, $0$ is the most uncomfort situation, while $1$ is the most comfort one. The $C^d_{j}$ is the desired water temperature in the $j$-th EWH. $\delta_{j}$ is the maximum allowed temperature deviation from the desired water temperature. The $C_{j}^{\text{min}}$ denotes the minimum temperature that can be tolerated by the $j$-th EWH.

\subsubsection{Electrical Storage}
In our proposed system, ES units can satisfy part of the power demand instantly through their discharging process. We use $E_{n, t}^s$ to denote the energy stored in the $n$-th ES unit at the end of period $t$, with the initial available energy $E_{n, 0}^s$. Then, we have the following dynamics for energy stored in the $n$-th ES unit:
\begin{equation}\label{con:es_dynamics1}
E_{n, t}^s = E_{n, t-1}^s + p_{n, t}^{s, +} \eta_{n}^{+} t - p_{n, t}^{s, -}/\eta_{n}^{-} t, \forall n, t, s,
\end{equation}
where $p_{n, t}^{s, +}$ and $p_{n, t}^{s, -}$ are the power charged into or discharged from the $n$-th ES unit at time $t$ in scenario $s$, and $\eta_{n}^{+}$ and $\eta_{n}^{-}$ represent the charging and discharging efficiencies of the $n$-th ES unit, respectively.

Each ES unit has a finite capacity, therefore, the energy stored in it must have the following lower and upper bounds:
\begin{equation}\label{con:es_dynamics2}
\underline{E}_n \leq E_{n, t}^s \leq \overline{E}_{n}, E_{n, 1}^s = E_{n, T}^s, \forall n, t, s,
\end{equation}
where $\overline{E}_{n}$ is the upper bound and $\underline{E}_n$ is the lower bound of the $n$-th ES unit's storage capacity. Moreover, we set the initial available energy the same as the final available energy for better scheduling in each operating day.

Furthermore, ES units have the charging and discharging rate limits as follows:
\begin{align}\label{con:es_dynamics3}
& 0 \leq p_{n, t}^{s, +} \leq P_{n}^{+} u_{n, t}^{s, +}, 0 \leq p_{n, t}^{s, -} \leq P_{n}^{-} u_{n, t}^{s, -}, \forall n, t, s \\
& 0 \leq u_{n, t}^{s, +} + u_{n, t}^{s, -} \leq 1, \forall n, t, s,
\end{align}
where $P_n^{+}$ and $P_n^{-}$ denote the maximum charged and discharged energy over period $t$ for the $n$-th ES unit, respectively. The $u_{n, t}^{s, +}$ and $u_{n, t}^{s, -}$ are binary variables indicating the charging and discharging decisions of the $n$-th ES unit. They are mutually exclusive, as given in~\eqref{con:es_dynamics3}.
\subsubsection{Heat Boiler}
Gas-based heat boilers can provide part of hot water demand instantly, while the other hot water demand is satisfied by EWHs~\citep{AhCh17}. Additionally, the utilization of heat boiler depends on the gas prices, and day-ahead and real-time retail electricity prices. Heat boilers generate the hot water either when the electricity prices are too high or when the heat balance cannot be maintained by EWHs alone. Thus, we have the following heat boiler generation equation:
\begin{equation}\label{eq:Heatoutput}
0 \leq h_{j, t}^s \leq \overline{H}_j, \forall j, t, s,
\end{equation}
where $h_{j, t}^s$ denotes the heat output of heat boiler at time $t$. The $\overline{H}_j$ is the upper bound of the heat output.

\subsection{Energy Balance}
We denote the aggregate critical power loads as $D_t^s$ that must be satisfied at each period $t$. Additionally, $N_i$, $N_i^k$, $N_i^j$, $N_i^n$, and $N_i^m$ represent sets of HVAC systems, PEVs, EWHs, ES units, and RES units of CB $i$, respectively.

\subsubsection{Power Balance}
The power balance equation is enforced in~\eqref{con:p_balance} as follows:
\begin{align}\label{con:p_balance}
& \sum_{i = 1}^{N_i} \sum_{n = 1}^{N_i^{n}} \left(p^{s, -}_{n, t} - p^{s, +}_{n, t}\right) + \sum_{i = 1}^{N_i} \sum_{m=1}^{N_i^m} w_{m, t}^s + g_{\text{RT}, t}^{s, +} \\
& = g_{\text{RT}, t}^{s, -} + D_t^{s} + \sum_{i = 1}^{N_i} \sum_{j = 1}^{N_i^j} l_{j, t}^s + \sum_{i = 1}^{N_i} P_{i, t}^{s, \text{hvac}} + \sum_{i = 1}^{N_i} \sum_{k = 1}^{N_i^{k}} p^{s, +}_{k, t}, \forall t, s. \nonumber
\end{align}

\subsubsection{Hot Water Balance}
As aforementioned, the hot water demand is satisfied by both EWHs and gas-based heat boilers, therefore we formulate the hot water balance equation as follows:
\begin{equation}\label{con:h_balance}
\zeta_j \sum_{i = 1}^{N_i} \sum_{j = 1}^{N_i^j} l_{j, t}^s + \sum_{i = 1}^{N_i} \sum_{j = 1}^{N_i^j} h_{j, t}^s \geq Q_t^s, \forall t, s,
\end{equation}
where $Q_t^s$ is the aggregated critical hot water demand. We also assume that the surplus heat can be disposed without penalty.

\subsection{Retail Electricity Market}
In our model, the central controller of CBs can either buy power from or sell power to a retail electricity market. Prior to an operating day, surplus power can be sold in a retail electricity market with the uncertain day-ahead selling price $c_{\text{DA}, t}^{s, -}$, and shortage power can be purchased from a retail electricity market with the uncertain day-ahead purchase price $c_{\text{DA}, t}^{s, +}$ at each time period $t$. We further denote $g_{\text{DA}, t}^{-}$ and $g_{\text{DA}, t}^{+}$ as the amount of power sold into and bought from a day-ahead retail electricity market at time $t$, respectively. After all uncertainties are revealed in an operating day, the central controller of CBs exchanges electricity in an operating day based on the real-time buying/selling electricity prices ( $c_{\text{RT}, t}^{s}$), where we assume the real-time buying and selling electricity prices are the same~\citep{LiCh17}. A penalty $\phi_t^s$ will occur when the electricity amount of the day-ahead schedule and the real-time delivery does not match~\citep{LiGu17}. We use $\overline{g}_t$ to represent the capacity limits on the point of common coupling (PCC). Moreover, in the business hours, a DSO may inform the central controller of CBs about reducing certain amount of power delivery through the PCC (in order to protect substations and transformers). Thus, $v_{t, s} \in [0.8, 1]$ is adopted to model the uncertain DR signals from a DSO in the buisness hours. Then, we have the following constraints for both day-ahead scheduled electricity amount and real-time delivery electricity amount:
\begin{align}\label{eq:connection1}
& 0 \leq g_{\text{DA}, t}^{-} \leq \overline{g}, 0 \leq g_{\text{DA}, t}^{+} \leq \overline{g}, \forall t \nonumber \\
& 0 \leq g_{\text{RT}, t}^{s, -} \leq \overline{g} v_t^s, 0 \leq g_{\text{RT}, t}^{s, +} \leq \overline{g} v_t^s, \forall t, s.
\end{align}

\section{Problem Formulation}\label{sec:problem}

\subsection{Two-Stage Stochastic Formulation}
\begin{figure}[!ht]
\centering
\includegraphics[width=4.8in]{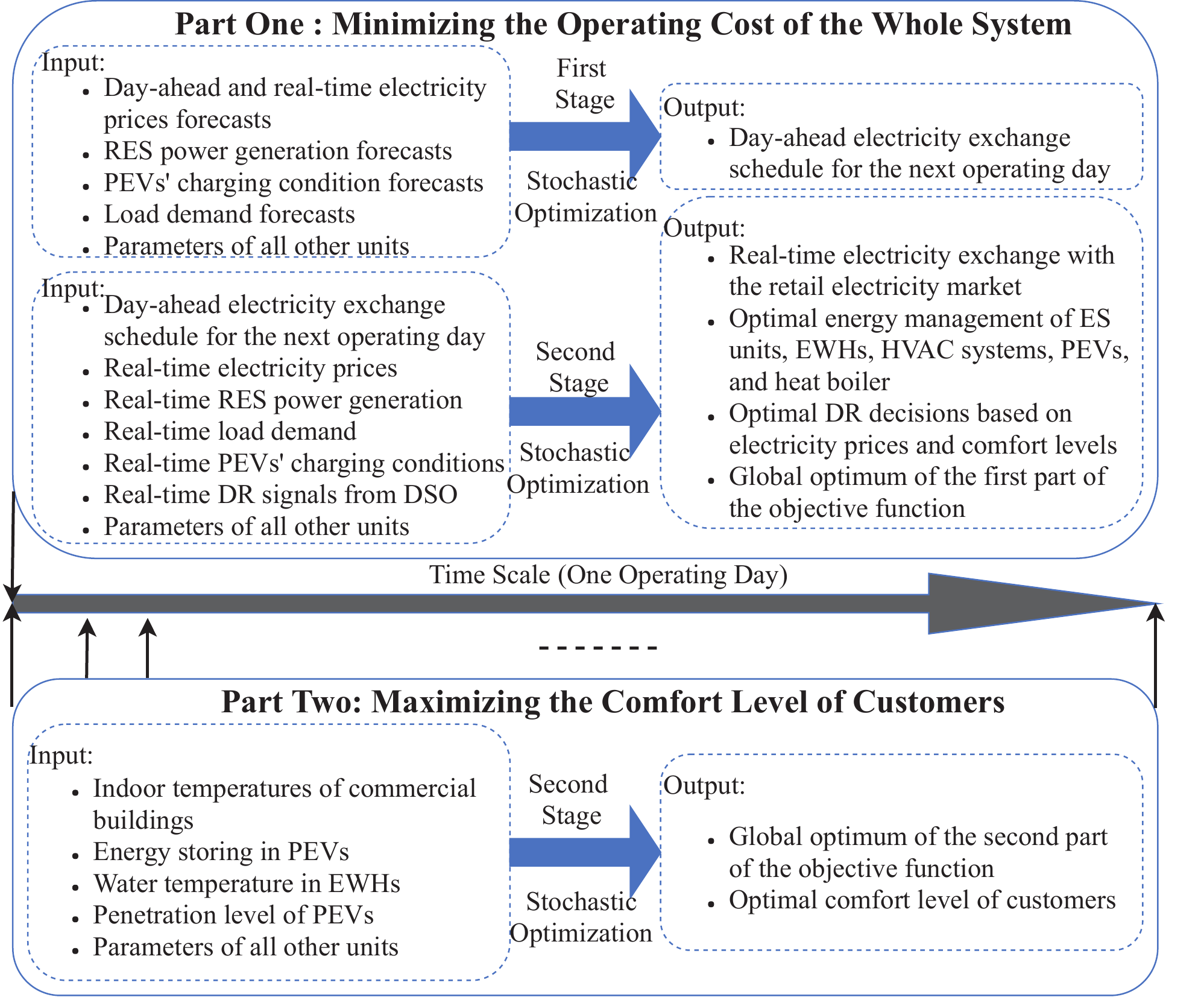}
\caption{A detailed two-stage stochastic programming.}
\label{fig:BilevelTwoStage}
\end{figure}
The objective of the proposed framework is to minimize the expected operating costs, while maximizing the comfort levels of consumers, which is shown in Fig.~\ref{fig:BilevelTwoStage}. This unique multi-objective problem can be formulated through a two-stage programming, with the following benefits: (i) our model is not only price-driven but also comfort-driven; and (ii) a two-stage modeling is more consistent with a two-settlement pool-based retail electricity market, where the unique features of the day-ahead electricity market and the real-time electricity market are reserved instead of simplifying any market modeling~\citep{MaJo18}.

\subsubsection{Objective Function Related with Operating Costs}
The first part of our two-stage stochastic formulation aims to minimize the expected operating costs of a commercial campus:
\begin{align}\label{eq:twostage}
& \min \sum_{s=1}^{N_s} \rho_s \sum_{t=1}^{N_t} \bigg\{ c_{\text{DA}, t}^{s, +} g_{\text{DA}, t}^{+} - c_{\text{DA}, t}^{s, -} g_{\text{DA}, t}^{-} \nonumber \\
& + c_{\text{RT}, t}^{s} \left[ g_{\text{RT}, t}^{s, +} - g_{\text{DA}, t}^{+} - \left(g_{\text{RT}, t}^{s, -} - g_{\text{DA}, t}^{-}\right) \right] \nonumber \\
& + c_t^p \left( \left| g_{\text{DA}, t}^{+} - g_{\text{RT}, t}^{s, +} \right| + \left| g_{\text{DA}, t}^{-} - g_{\text{RT}, t}^{s, -} \right| \right) \nonumber \\
& + \sum_{i = 1}^{N_i} \sum_{k=1}^{N_i^{k}} c_k p_{k, t}^{s, +} \eta_k^{+} \nonumber \\
& + \sum_{i = 1}^{N_i} \sum_{n=1}^{N_i^{n}} c_n \left(p_{n, t}^{s, +} \eta_n^{+} + p_{n, t}^{s, -}/\eta_n^{-} \right) \nonumber \\
& + \sum_{i = 1}^{N_i} \sum_{j = 1}^{N_i^j} c_t^b h_{j, t}^s \bigg \},
\end{align}
subject to constraints~\eqref{con:HVAC1}--\eqref{con:HVAC2}, \eqref{con:pev_dynamics1}--\eqref{con:pev_dynamics3}, \eqref{con:EWH1}--\eqref{con:EWH4}, and \eqref{con:es_dynamics1}--\eqref{eq:connection1}.

The first line represents the day-ahead bidding cost of the central controller of CBs, the second line shows the real-time electricity cost of the central controller of CBs, the third line is the penalty cost for differences between the day-ahead and real-time electricity buy from or sell to a retail electricity market, the fourth line expresses the degradation costs of PEVs, the fifth line denotes the degradation costs of ES units, and the last line is the generation costs of gas-based heat boilers.

\subsubsection{Objective Function Related with Comfort Levels}
The second part of our two-stage stochastic formulation focuses on maximizing the expected comfort levels related to HVAC systems, PEVs and EWHs. The comfort levels for HVAC systems and EWHs consider the overall population of each CB, while the comfort levels for PEVs considers only the PEV owners. We define the penetration levels of PEVs to be $PL_{pev}=N_i^k/N_i^p$, where $N_i^p$ is the total population in the $i$-th CB.
\begin{equation}\label{eq:twostage2}
\max \sum_{s=1}^{N_s} \rho_s \sum_{t=1}^{N_t} \sum_{i=1}^{N_i} \left( J_{i, t}^s + \frac{N_i^k}{N_i^p} \sum_{k=1}^{N_i^{k}} J_{k, t}^s + \sum_{j=1}^{N_i^j} J_{j,t}^s \right)
\end{equation}
subject to constraints~\eqref{con:comfort1},~\eqref{con:comfort2} and \eqref{con:comfort3}.

\section{Solution Methodology}\label{sec:solution}

\subsection{Problem Reformulation}
To integrate these two objectives, we reformulate the first objective as maximizing the expected benefits of CBs. Then we add it to the second objective after unification. An unification process is coupled with comfort levels, where the total power that is required to increase all comfort levels from 0 to 1, is calculated as base power $P_{base}$. Therefore, we have the following reformulated two-stage stochastic programming:
\begin{align}\label{eq:reformulation}
& \max \sum_{s=1}^{N_s} \rho_s \sum_{t=1}^{N_t} \bigg\{ c_{\text{DA}, t}^{s, -} g_{\text{DA}, t}^{-}/P_{base} - c_{\text{DA}, t}^{s, +} g_{\text{DA}, t}^{+}/P_{base} \nonumber \\
& + c_{\text{RT}, t}^{s} \bigg\{ \left(g_{\text{RT}, t}^{s, -}/P_{base} - g_{\text{DA}, t}^{-}/P_{base} \right) - g_{\text{RT}, t}^{s, +}/P_{base} \nonumber \\
& + g_{\text{DA}, t}^{+}/P_{base}\bigg\} - c_t^p \psi_t^s \nonumber \\
& - \sum_{i = 1}^{N_i} \sum_{k=1}^{N_i^{k}} c_k p_{k, t}^{s, +} \eta_k^{+} /P_{base} \nonumber \\
& - \sum_{i = 1}^{N_i} \sum_{n=1}^{N_i^{n}} c_n \left(p_{n, t}^{s, +} \eta_n^{+} /P_{base} + p_{n, t}^{s, -}/\eta_n^{-} /P_{base} \right) \nonumber \\
& - \sum_{i = 1}^{N_i} \sum_{j = 1}^{N_i^j} c_t^b h_{j, t}^s/\overline{H}_j \nonumber \\
& + \sum_{i=1}^{N_i} J_{i, t}^s + \sum_{i = 1}^{N_i} \frac{N_i^k}{N_i^p} \sum_{k=1}^{N_i^{k}} J_{k, t}^s + \sum_{i = 1}^{N_i} \sum_{j=1}^{N_i^j} J_{j,t}^s \bigg \},
\end{align}
subject to constraints~\eqref{con:HVAC1}--\eqref{eq:connection1}. $\psi_t$ is an auxiliary variable that replace the absolute term $\left( \left| g_{\text{DA}, t}^{+} - g_{\text{RT}, t}^{s, +} \right| + \left| g_{\text{DA}, t}^{-} - g_{\text{RT}, t}^{s, -} \right| \right)$ in the previous objective function~\eqref{eq:twostage}, where $\psi_t^s = \psi_{1,t}^s + \psi_{2,t}^s $, $\psi_{1,t}^s=\left| g_{\text{DA}, t}^{+} - g_{\text{RT}, t}^{s, +} \right|$ and $\psi_{2,t}^s=\left| g_{\text{DA}, t}^{-} - g_{\text{RT}, t}^{s, -} \right|$.

\subsection{Stochastic Programming}
We use a scenario-based two-stage stochastic programming to handle uncertainties related to the output of renewables, PEV charging conditions, uncertain load demand, business hour DR signal from the DSO and day-ahead and real-time electricity prices. Uncertainties are modeled by scenarios, where the historic data of various uncertainties can be obtained by the central controller of CBs. Also, the day-ahead electricity prices and real-time electricity prices are correlated in real-world practice, therefore, scenarios regarding the day-ahead and real-time electricity prices are generated together. Other uncertainties are assumed to be independent of each other. Then, we use the historic data of previous 30 days to generate $30 \times 30 \times 30 \times 30 = 810,000$ scenarios with the same probability through Monte Carlo method. However, this huge number of scenarios require high computational efforts and make this optimization problem intractable. Therefore, a scenario reduction method is necessary to reduce the original amount of scenarios to a reasonable quantity~\citep{WaLi17}.

\subsection{Scenario Reduction and Linearization}
To achieve this goal, we use a fast-forward scenario reduction method to reduce the original $810,000$ scenarios to 30 scenarios. Our objective is to obtain a small set of scenarios that can maintain the properties of original scenarios at a reasonable level.

As observed from the optimization model, there exists two nonlinear terms in the objective function: $\psi_{1,t}^s=\left| g_{\text{DA}, t}^{+} - g_{\text{RT}, t}^{s, +} \right|$ and $\psi_{2,t}^s=\left| g_{\text{DA}, t}^{-} - g_{\text{RT}, t}^{s, -} \right|$. In order to linearize the aforementioned nonlinear terms, we adopt the following constraints:
\begin{align}
& \psi_{1,t}^s \geq g_{\text{DA}, t}^{+} - g_{\text{RT}, t}^{s, +}, \psi_{1,t}^s \geq g_{\text{RT}, t}^{s, +} - g_{\text{DA}, t}^{+}, \forall t, s \\
& \psi_{2,t}^s \geq g_{\text{DA}, t}^{-} - g_{\text{RT}, t}^{s, -}, \psi_{2,t}^s \geq g_{\text{RT}, t}^{s, -} - g_{\text{DA}, t}^{-}, \forall t, s.
\end{align}
Therefore, this two-stage stochastic programming can be formulated as a large-scale mixed integer linear program (MILP), which can be directly solved by commercial solvers such as GUROBI. Specifically, there are  $37,440$ equality constraints and $67,680$ inequality constraints in the optimization model, with $207,015$ continuous variables and $17,280$ integer variables ($17,280$ binary variables).

\section{Simulation Results}\label{sec:simulation}
In this section, the proposed algorithm is evaluated through the real-world datasets. Firstly, the implemented real-world data sets are described in detail. Then, the performances of the proposed stochastic solution are evaluated.

All simulations are implemented on a desktop computer with 3.0 GHz Intel Core i5-7400 CPU and 8GB RAM. All the scenario generation and scenario reduction processes are performed using the MATLAB. The proposed energy management problem are simulated with the python 2.7. A convergence criterion is set as the MIPGAP default value ($0.0001$) in GUROBI 8.0.0.

\subsection{Numerical Settings}
\begin{table}[!ht]
\centering
\caption{Comfort Level Related Parameters}
\label{table:comfortlevel}
\begin{tabular}{|c| |c| |c| |c|}
\hline
Type & $C^d_{i}$ ($^\circ C$) & $\delta_{i}$ ($^\circ C$) & $\epsilon_{i}$ ($^\circ C$) \\
\hline
HVAC & 24 & 2 & 0.5  \\
\hline
Type & $E_{k}^d$ ($\%$) & $E_{k}^{\text{base}}$ ($\%$) & $E_k^0$ ($\%$) \\
\hline
PEV & 80 & 10 & 10 \\
\hline
Type & $C_{j}^d$ ($^\circ C$) & $\delta_{j}$ ($^\circ C$) & $C_{j}^{0}$ ($^\circ C$)\\
\hline
EWH & 40 & 10 & 30 \\
\hline
\end{tabular}
\end{table}
\begin{table}[!ht]
\centering
\caption{HVAC System Parameters}
\label{table:CB}
\begin{tabular}{|c| |c| |c| |c| |c| |c| |c|}
\hline
HVAC & CB1  & CB2 &  CB3 & CB4 & CB5 & CB6\\
\hline
$\overline{P}_{i}^{\text{hvac}}$ & 0.1 & 0.15 & 0.2 & 0.25 & 0.3 & 0.35  \\
\hline
\end{tabular}
\end{table}
\begin{table}[!ht]
\centering
\caption{ES Unit Parameters}
\label{table:3}
\begin{tabular}{|c| |c| |c| |c|}
\hline
$\underline{E}_n$ (kWh) & $\overline{E}_{n}$ (kWh) & $P_{n}^{+}$ (kW) & $P_n^{-}$ (kW) \\
\hline
4 & 76 & 4 & 4    \\
\hline
\end{tabular}
\end{table}
\begin{table}[!ht]
\centering
\caption{PEV Parameters}
\label{table:pev}
\begin{tabular}{|c| |c| |c| |c|}
\hline
Type & $\underline{E}_k$ (kWh) & $\overline{E}_{k}$ (kWh) & $P_{k}^{+}$ (kW) \\
\hline
Tesla Model S 75D & 3.8 & 75 & 11.5 \\
\hline
Tesla Model X 100D & 5 & 100 & 17.2 \\
\hline
Nissan Leaf SV & 1.5 & 30 & 3.6 \\
\hline
\end{tabular}
\end{table}
\begin{figure}[!ht]
\centering
\includegraphics[width=4.8in]{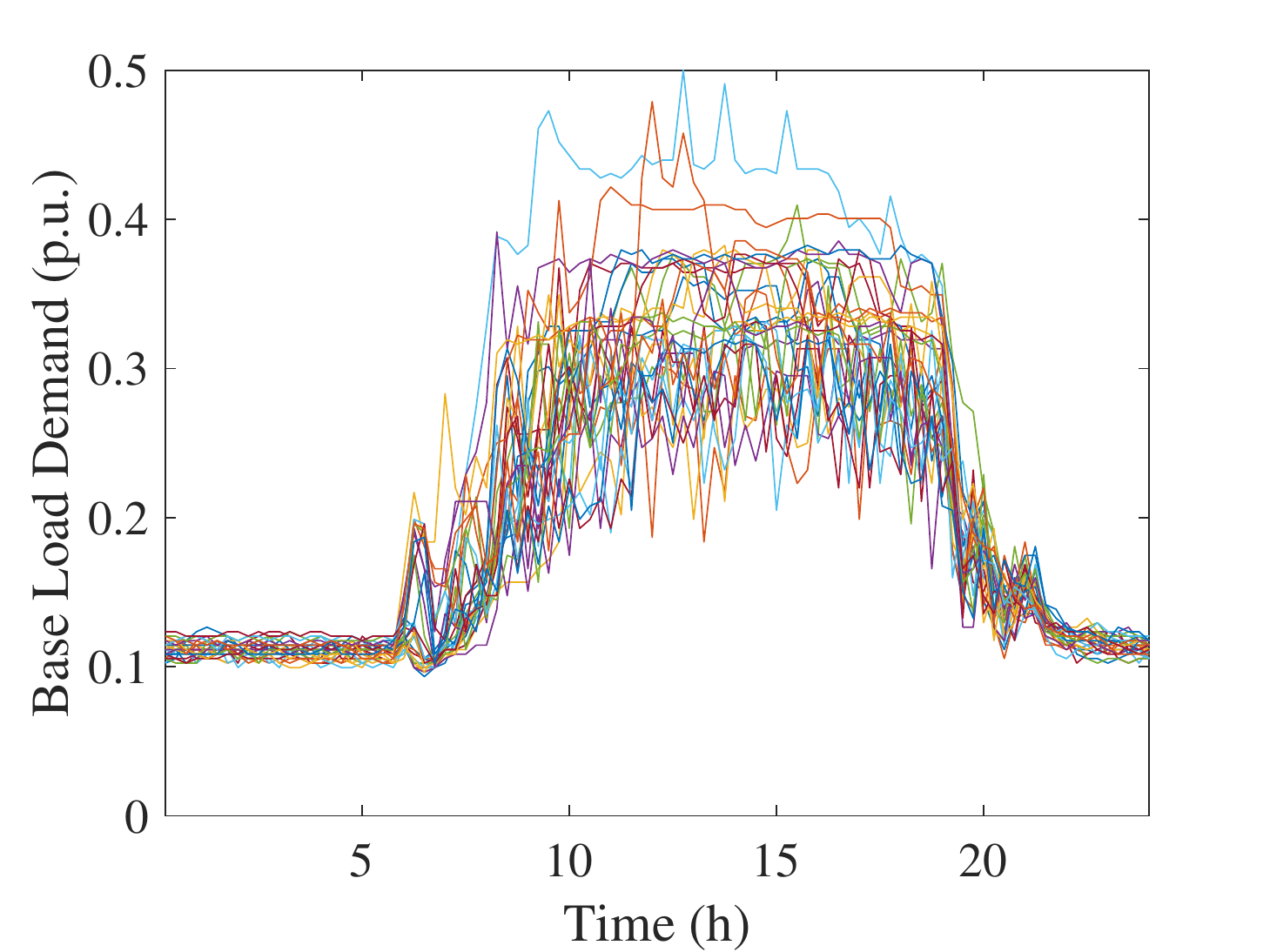}
\caption{30 scenarios of unified base load demand.}
\label{fig:load}
\end{figure}
\begin{figure}[!ht]
\centering
\includegraphics[width=4.8in]{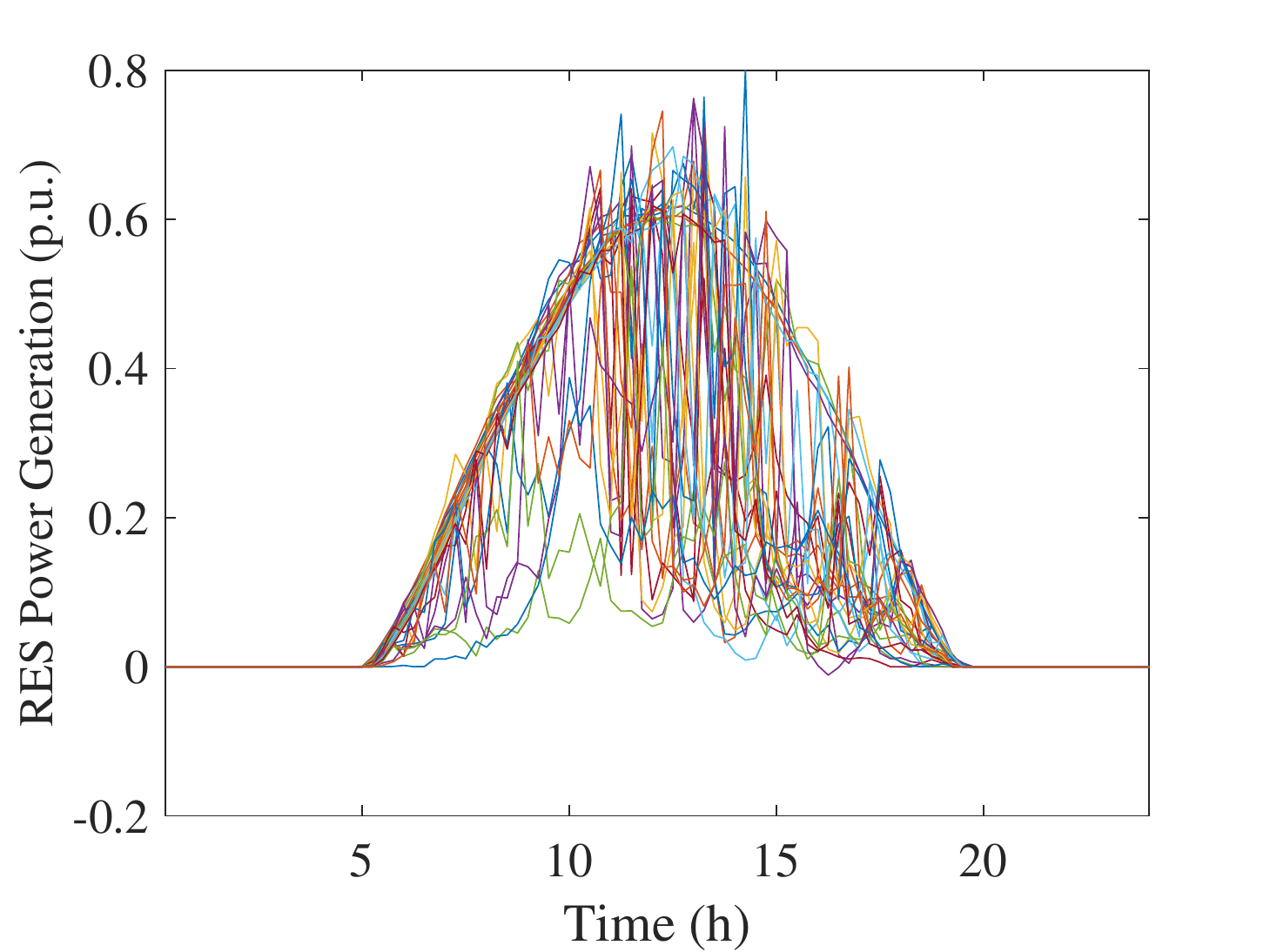}
\caption{30 scenarios of unified output of solar panels.}
\label{fig:solar}
\end{figure}
The considered commercial campus consists of six CBs, six packs of roof-top solar panels, fifty PEVs in three different types, one aggregated critical power load, and one aggregated base hot water load. Each CB includes one EWH, one HVAC system , one pack of batteries and one gas-based heat boiler, which is scalable. All parameters are unified for the computational proposes, with $P_{\text{base}}$ as $1867 kW$ and $H_{\text{base}}$ as $1224 kBtu$. Comfort levels are only considered during the business hours when consumers are on the commercial campus (from $8$ a.m. to $8$ p.m.), where the parameters are given in Table~\ref{table:comfortlevel}.

Among all the real-world data, data for the base power and hot water loads are taken from the Global Energy Interconnection Research Institute North America (GEIRINA)'s CBs in the San Jose campus, as shown in Fig.~\ref{fig:load}. Each heat boiler has a maximum generation of $206 kBtu$, where the gas price is also from the GEIRINA's CBs in the San Jose campus. The total installed capacity of the solar panels is $1500 kW$, where the historic generation patterns are taken from~\citep{Wind}, as shown in Fig.~\ref{fig:solar}. Also, the historical data of the solar irradiance and the outdoor temperature are from~\citep{Wind}, with proper scaling coefficients. The day-ahead and real-time electricity prices are taken from the real-world wholesale electricity prices of PJM~\citep{PJM15}, with proper scaling coefficients.

As shown in Table~\ref{table:CB}, each HVAC system has a different power upper bound. $\sigma_{i, t}$ is set to be $-1$, which makes HVAC systems perform cooling only. $\eta_i$, $\alpha_i$ and $\beta_i$ are taken from~\citep{NgLe14Comfort}. As shown in Table~\ref{table:3}, each battery pack has a storage capacity of $80 kWh$, and has charging and discharging efficiencies of $0.98$. The maximum charging and discharging rates for each battery pack are both $4 kW$. To prolong the battery packs' lifetime, the energy levels should not drop below $5\%$ or overcharged above $95\%$ of the capacity. Both the initial and final stored electrical energy are set to be $50\%$ of its total capacity. Degradation cost coefficients $c_n$ and $c_k$ are set to be $0.0035 \$/kWh$. For PEVs, the major brands of PEVs are listed in Table~\ref{table:pev}, which share $70\%$, $10\%$ and $2\%$ of the total PEVs market in the state of California, respectively. The initial energy levels of all PEVs are set as $10\%$ of its total capacity. The day-ahead selling prices are set to be $80\%$ of the day-ahead purchasing prices. The electricity exchange limit between the proposed system and the retail electricity market is set to be $1867 kW$. The data for EWHs are taken from~\citep{SaOr15}. The EWHs' power to heat ratio is set to $1.2$. Each heat boiler has a maximum heat generation of $174 kBtu$.

\subsection{Case Study}

\subsubsection{First-Stage Decisions: With/Without the Comfort Level}
\begin{figure}[!ht]
\centering
\includegraphics[width=4.8in]{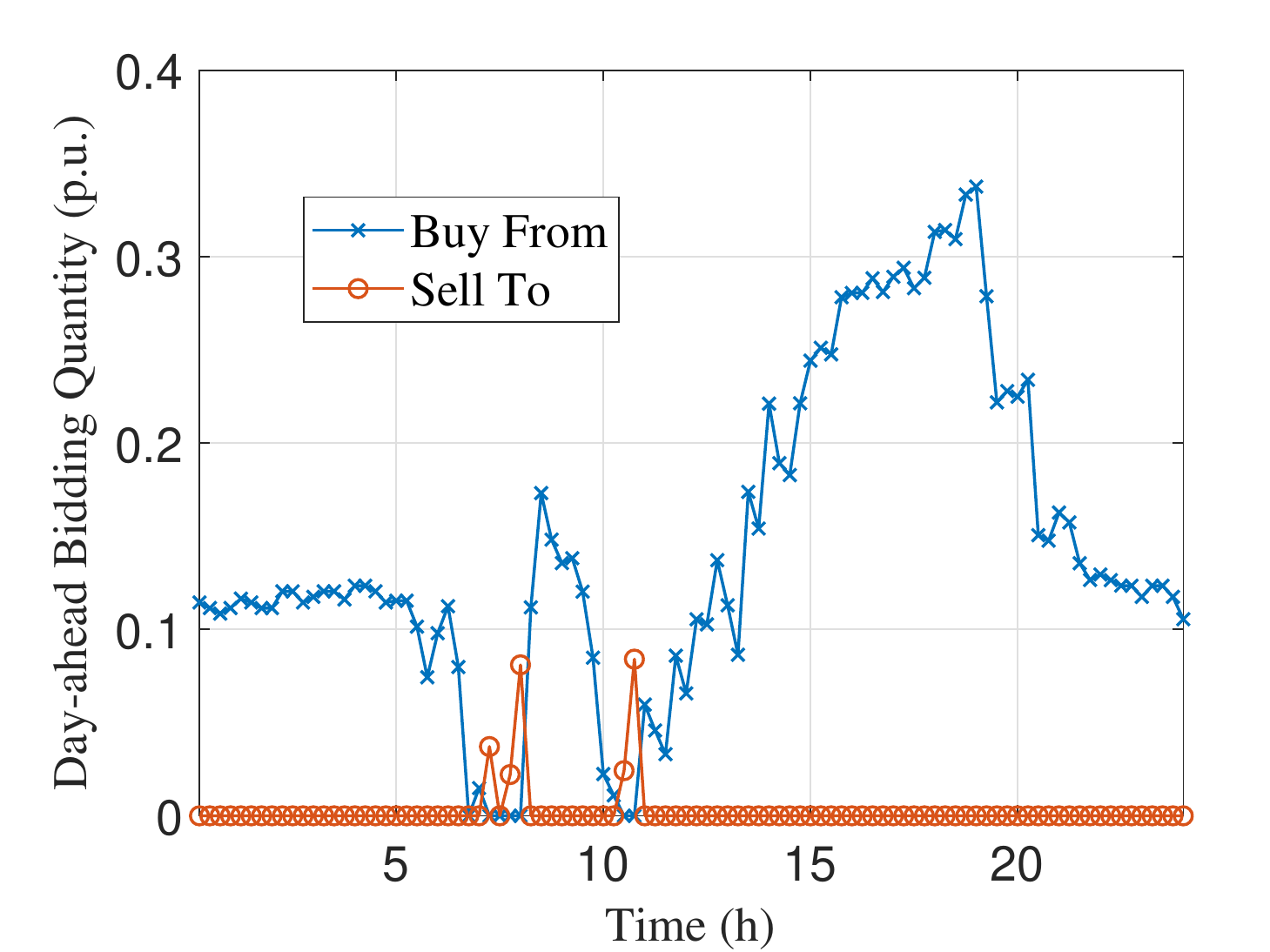}
\caption{Day-ahead electricity buy from or sell to the retail electricity market considering comfort levels.}
\label{fig:DA_withJ}
\end{figure}
\begin{figure}[!ht]
\centering
\includegraphics[width=4.8in]{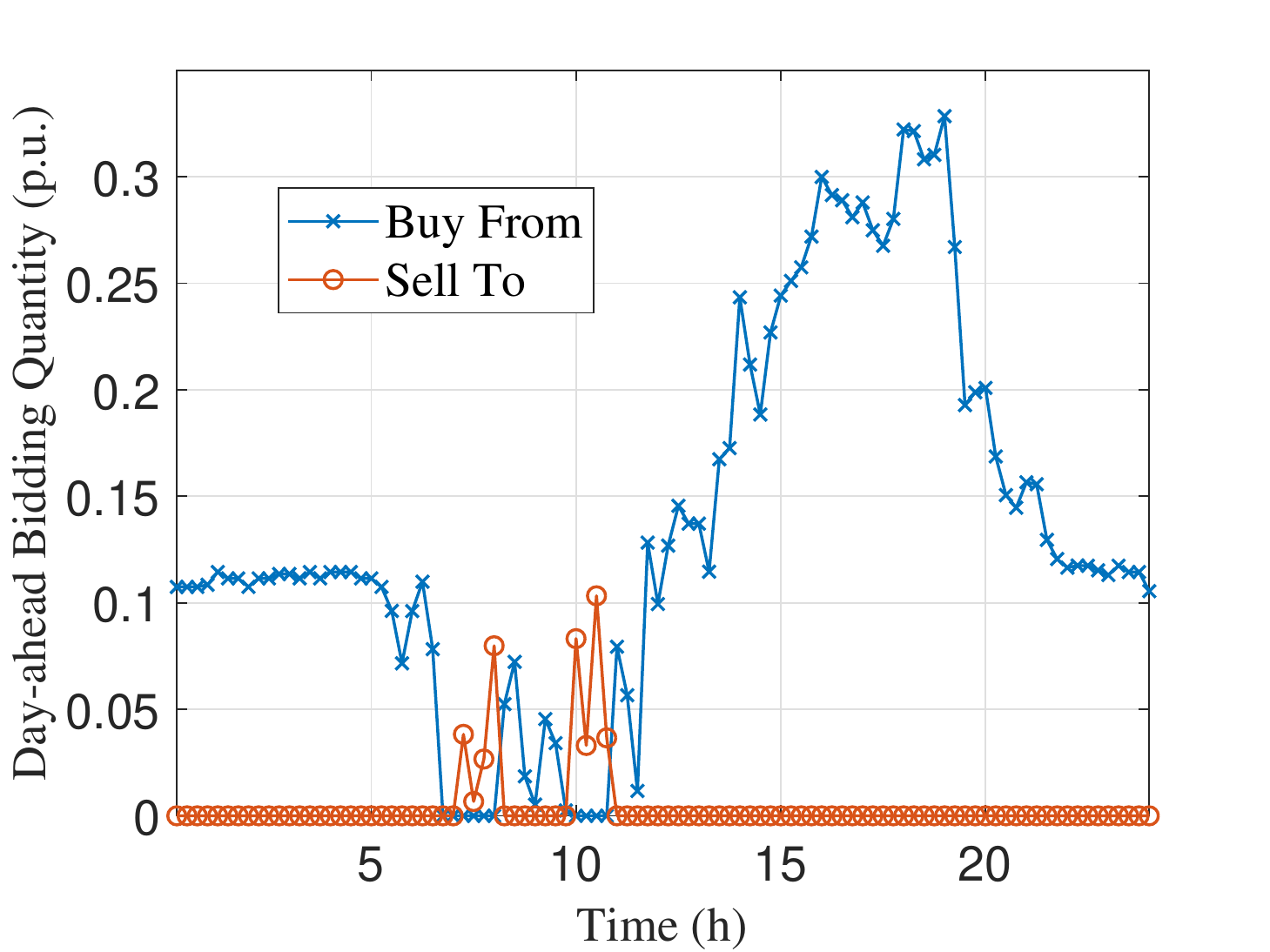}
\caption{Day-ahead electricity buy from or sell to the retail electricity market without comfort levels.}
\label{fig:DA_noJ}
\end{figure}
The proposed two-stage scenario-based stochastic optimization approach is tested based on two different cases, namely with/without the comfort level related objectives and constraints. As shown in Figs.~\ref{fig:DA_withJ} and~\ref{fig:DA_noJ}, the $x$-axis is the time (24 hours) of one operating day with $96$ time slots, while the $y$-axis shows unified numbers denoting the day-ahead electricity buy from or sell to the retail electricity market. An identical system model of the day-ahead bidding strategy is implemented for both cases (with and without the comfort level), thus, the base loads have the most significant influence on the energy usage when consumers are off-campus (from 0-8 and 20-24), which makes non-business hours biding pattern similar. During business hours when all consumers are on the campus (from 8 to 20), the comfort levels of HVAC systems, PEVs and EWHs play an important role on the central controller making decisions on the day-ahead bidding strategy.

It can also be observed that in case one, the central controller of CBs buy almost the maximum allowed power from the retail electricity market from $8$ to $11$ to support comfort level requirements, while selling a little amount of power to the retail electricity market when the solar panels can support sufficient energy to other appliances. However, for case two, the central controller of CBs only buy power when the electricity prices are low and sell a large amount of power when the solar panels have the maximum power outputs during business hours. This is reasonable since in case one, the central controller of CBs provide DR considering both electricity prices and comfort levels, which need to find a trade-off between minimizing operating costs and maximizing comfort levels; while in case two, the central controller of CBs are driven by only electricity prices, which ignores comfort needs of consumers.
\begin{table}[!ht]
\centering
\caption{Expected Total Operating Cost and Comfort Level of the Commercial Campus}
\label{table:4}
\begin{tabular}{|c| |c| |c| |c|}
\hline
& Without  & With  & $Increment\%$ \\
\hline
Expected Total Operating Cost & 11.45 & 13.44 & $17.4\%$\\
\hline
Expected Total Comfort Level & 0.4 & 1 & $150\%$\\
\hline
\end{tabular}
\end{table}
Interestingly, the pure unified expected operating costs of case one is $13.44$ compared with $11.45$ of case two, which is only $17.4\%$ of increase in costs. However, the expected comfort levels rise from $40\%$ to $100\%$ for all the appliances for the case considering comfort levels, as shown in Table~\ref{table:4}.

\subsubsection{Sensitivity Analysis: With/Without the Comfort Level}
\begin{figure}[!ht]
\centering
\includegraphics[width=4.8in]{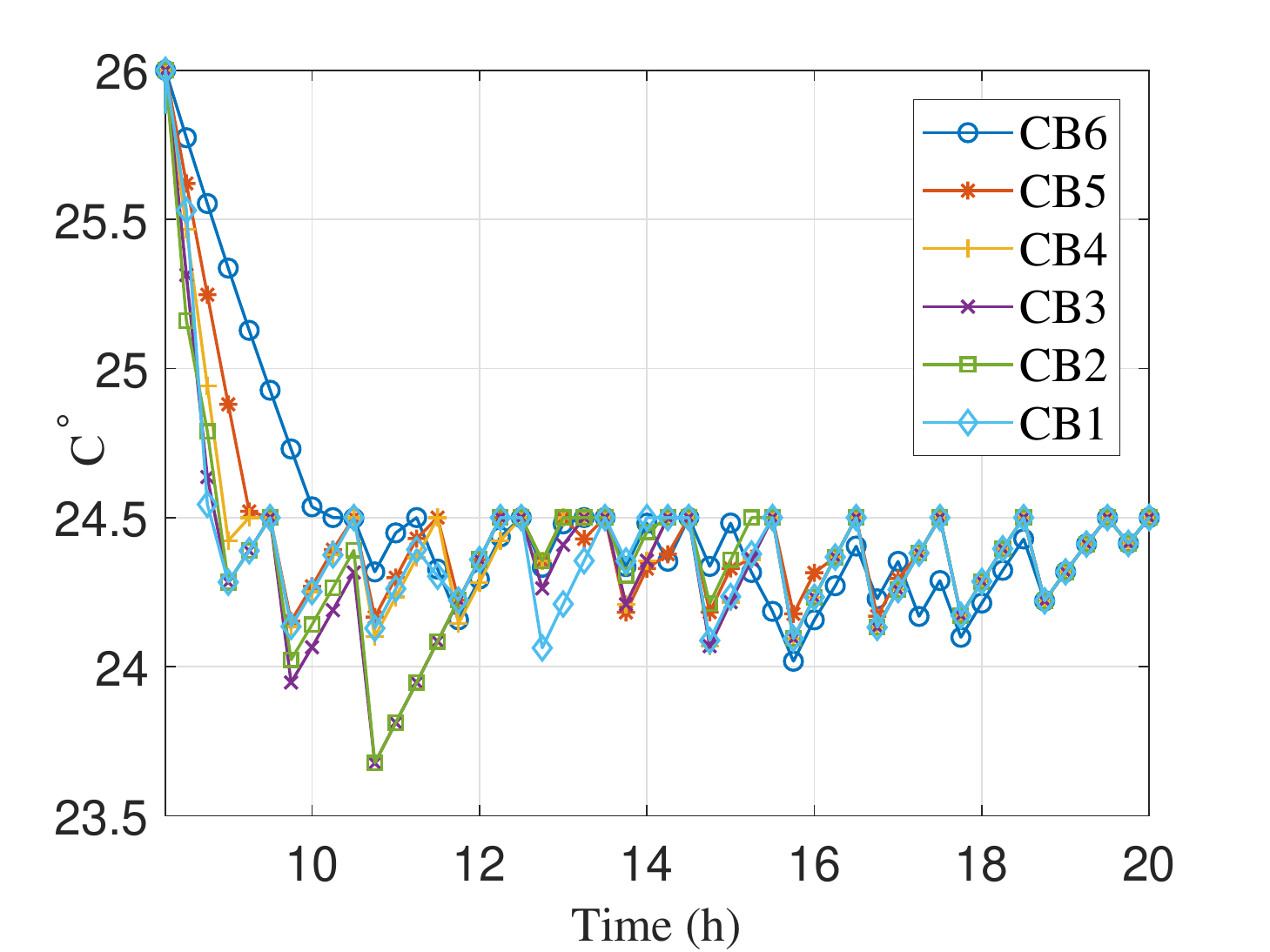}
\caption{Indoor temperatures with our proposed comfort level model.}
\label{fig:temp1}
\end{figure}
\begin{figure}[!ht]
\centering
\includegraphics[width=4.8in]{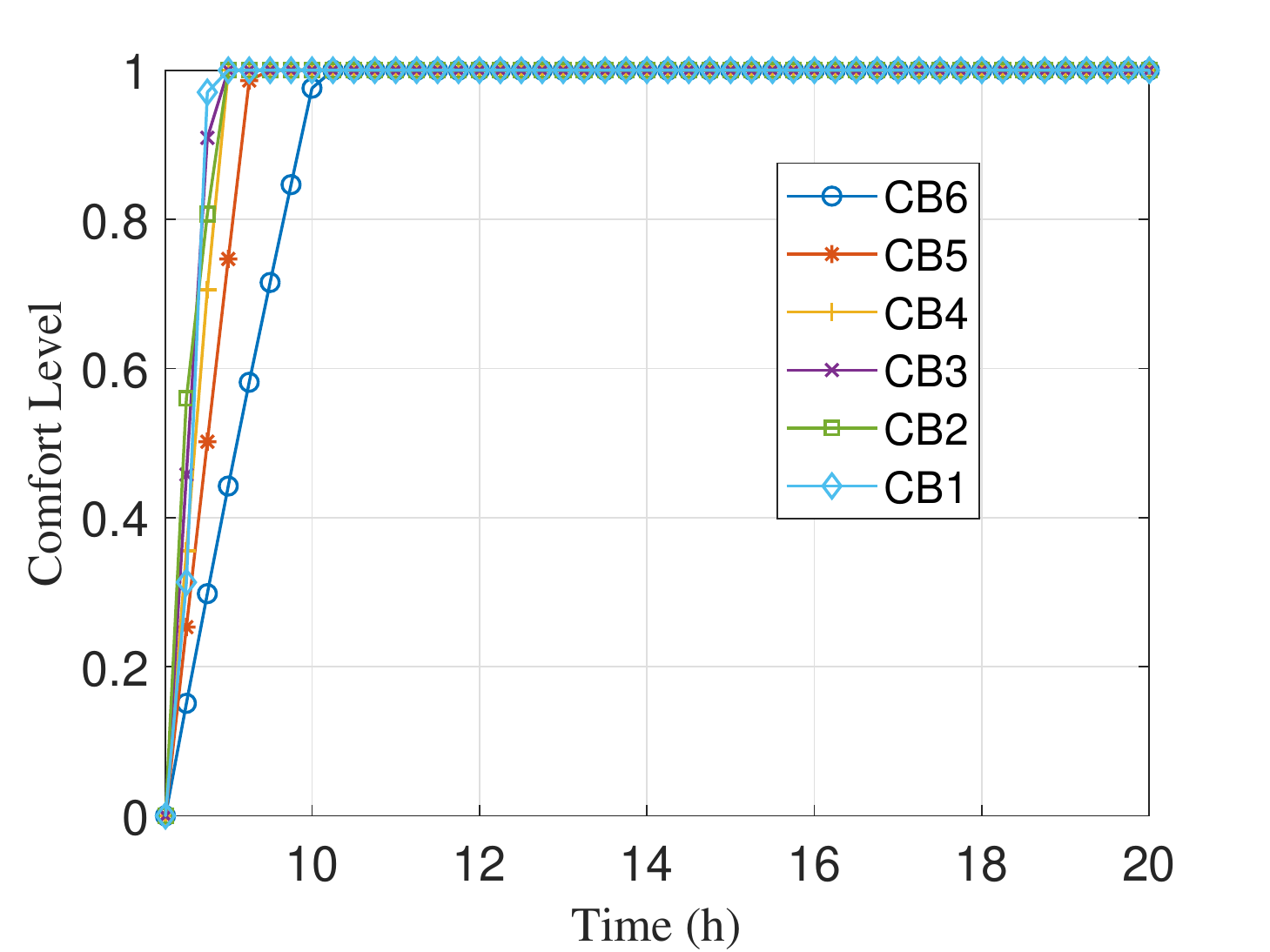}
\caption{Comfort levels related with HVAC systems with our proposed comfort level model.}
\label{fig:HVAC_D1}
\end{figure}
\begin{figure}[!ht]
\centering
\includegraphics[width=4.8in]{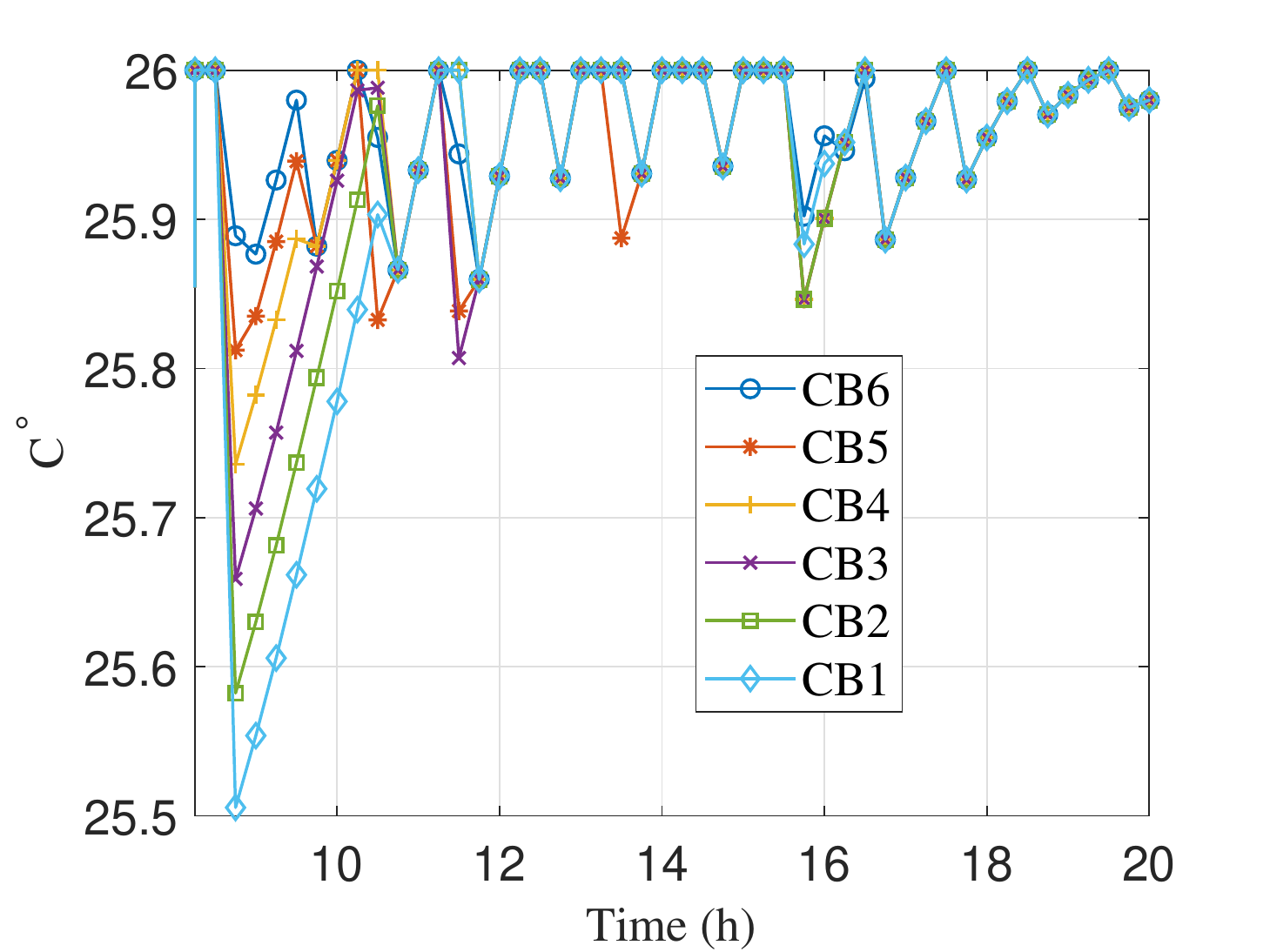}
\caption{Indoor temperatures without our proposed comfort level model.}
\label{fig:temp2}
\end{figure}
\begin{figure}[!ht]
\centering
\includegraphics[width=4.8in]{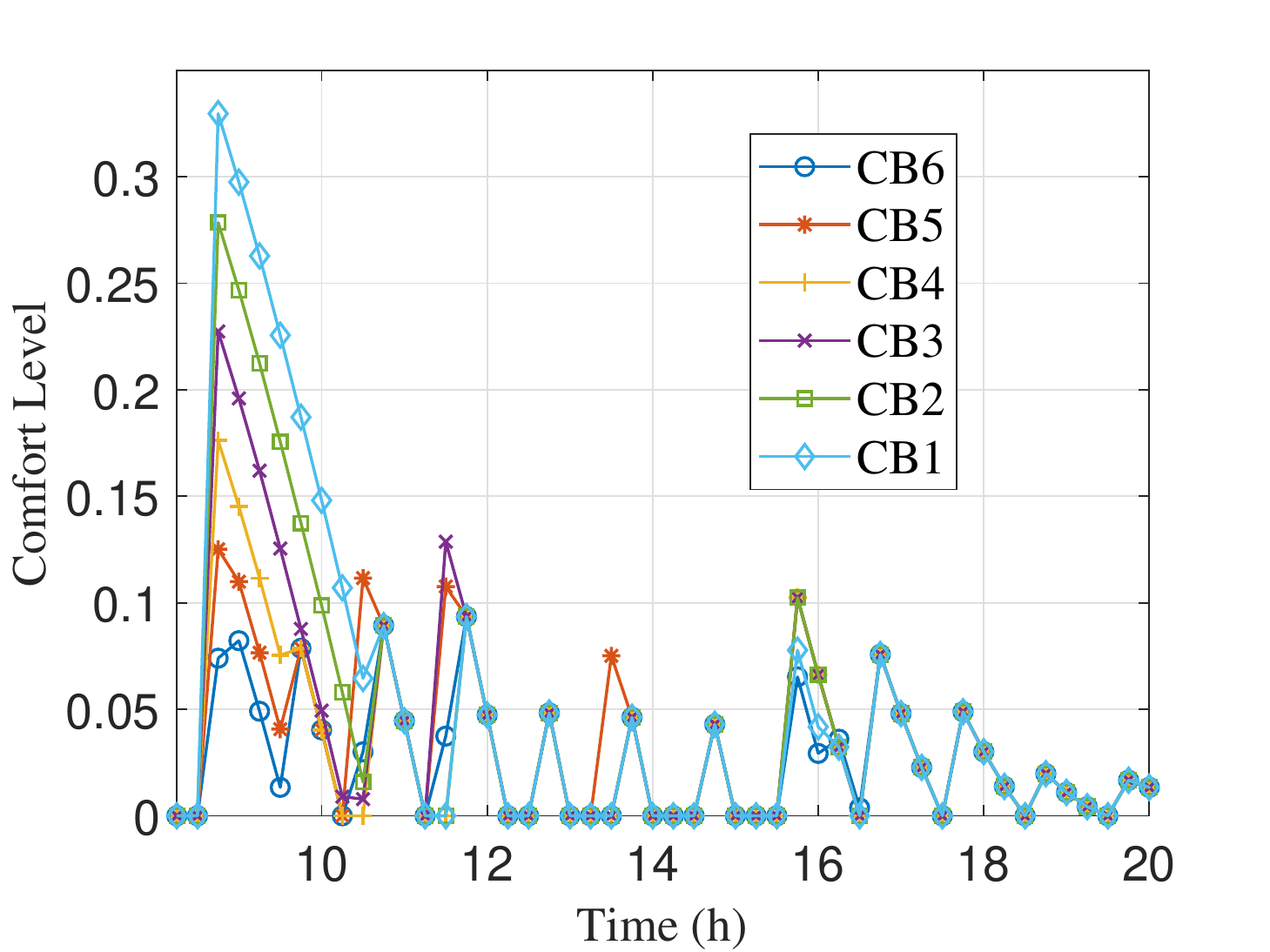}
\caption{Comfort levels related with HVAC systems without our proposed comfort level model.}
\label{fig:HVAC_D2}
\end{figure}
In order to show the effectiveness of the proposed comfort level models, the system is tested in a single scenario with and without comfort levels related to the objective functions, while comfort levels are still calculated through the original constraints. As shown in Figs.~\ref{fig:temp1} and~\ref{fig:temp2}, the indoor temperatures of six buildings reduce significantly when considering the proposed comfort level formulations, which is because the real-time electricity prices during this period (from 8 am to 10 am) are lower enough to provide sufficient power to HVAC systems to reduce indoor temperatures. This situation is also provided in Figs.~\ref{fig:HVAC_D1} and~\ref{fig:HVAC_D2}, where the comfort levels increase vastly during the same period. In the business hours when the solar irradiance is at its maximum radiation, the indoor temperatures increase to their upper bounds and comfort levels decrease accordingly for the other case. This is because even though the solar panels are also at their maximum outputs, however, the indoor temperatures are not the central controller's priority when the comfort levels are not considered. The indoor temperatures are maintained within the most comfort region with the help of the proposed HVAC comfort level models; however, for the one without comfort levels, the indoor temperatures remain high and the comfort levels keep low until sunset.

\begin{figure}[!ht]
\centering
\includegraphics[width=4.8in]{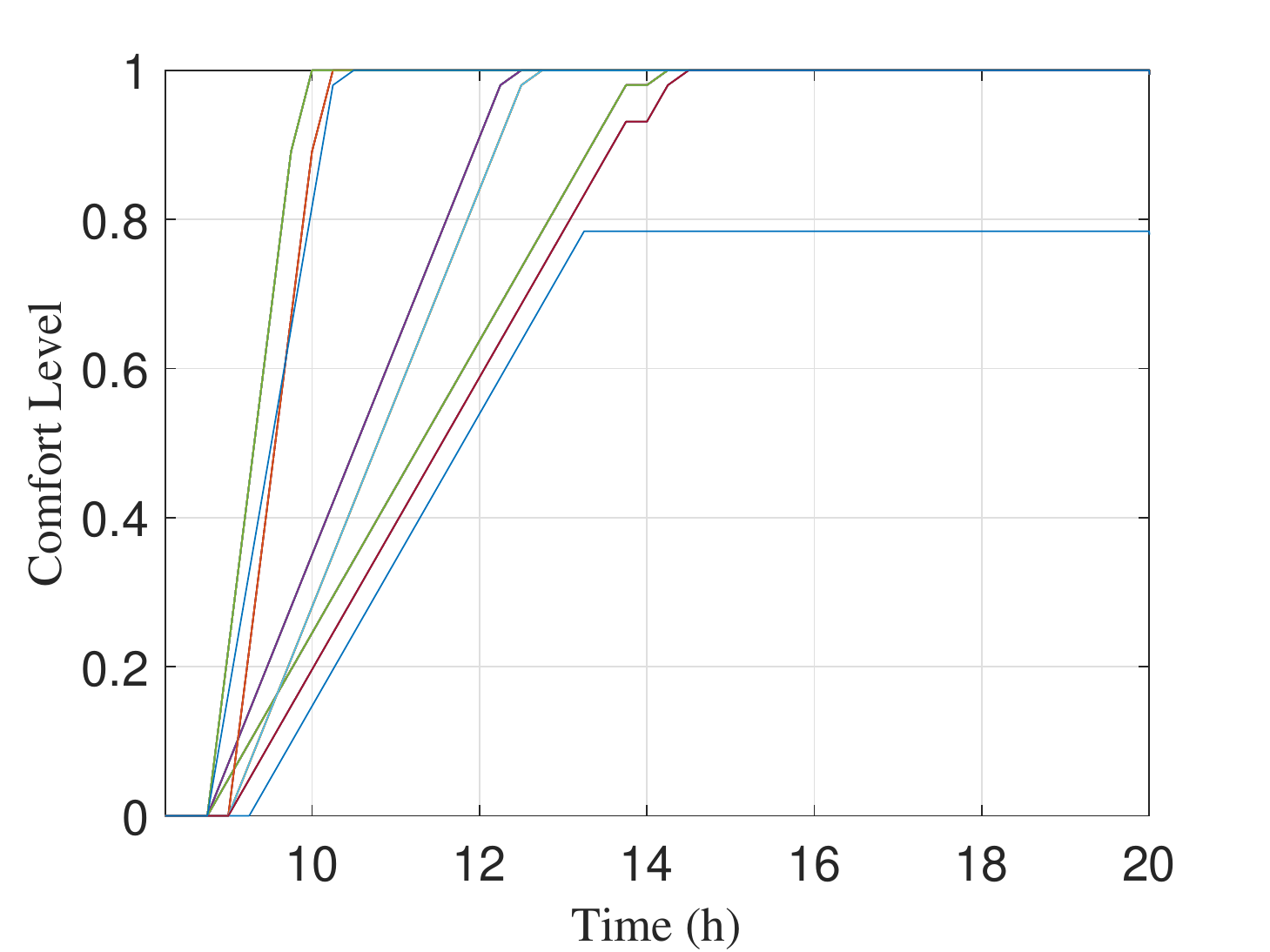}
\caption{Comfort levels related with PEVs considering our proposed comfort level model.}
\label{fig:PEV_D}
\end{figure}
\begin{figure}[!ht]
\centering
\includegraphics[width=4.8in]{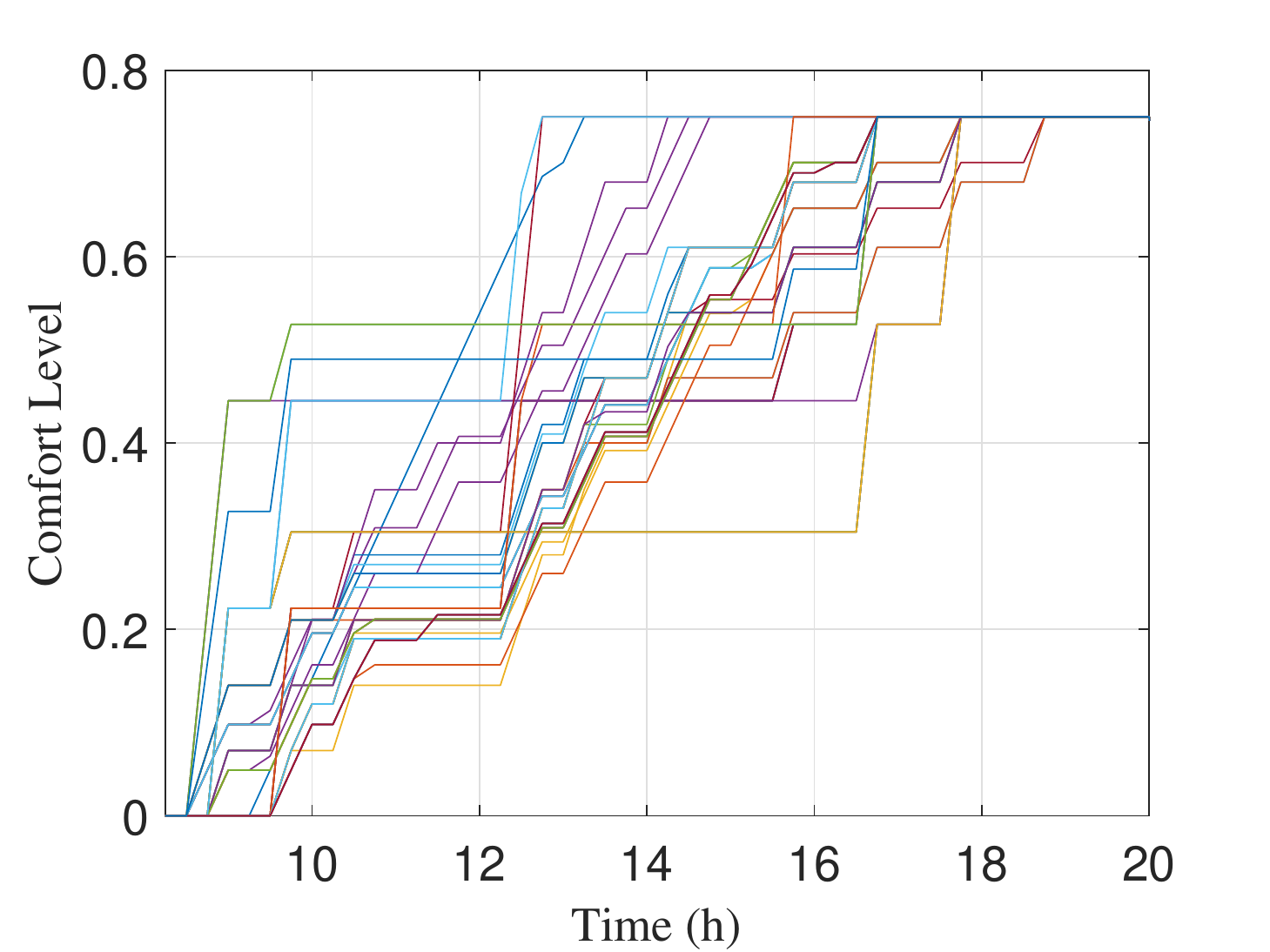}
\caption{Comfort levels related with PEVs without our proposed comfort level model.}
\label{fig:PEV_DnoJ}
\end{figure}
The comfort level models of PEVs and EWHs in the same circumstances with and without comfort levels are further tested. As shown in Figs.~\ref{fig:PEV_D} and~\ref{fig:PEV_DnoJ}, the comfort levels increase faster when considering the proposed comfort level models, however, this does not mean that the real-time electricity prices are not considered. Instead, both real-time electricity prices and comfort levels are considered to charge PEVs to the desired energy levels. In the other case, PEVs are charged to $80\%$ of the desired energy levels when the comfort level models are not in use. However, in this situation, only the electricity price impact is taken into consideration, which may result in less comfort when the consumers left earlier. As a result, comfort levels of these PEV owners are decreased. As shown in Fig.~\ref{fig:PEV_D}, the PEVs left at $13:15$ are charged to $97\%$ of their desired energy levels; while in the second case as shown in Fig.~\ref{fig:PEV_DnoJ}, these PEVs can only be charged to $42\%$ of their desired energy levels.

\begin{figure}[!ht]
\centering
\includegraphics[width=4.8in]{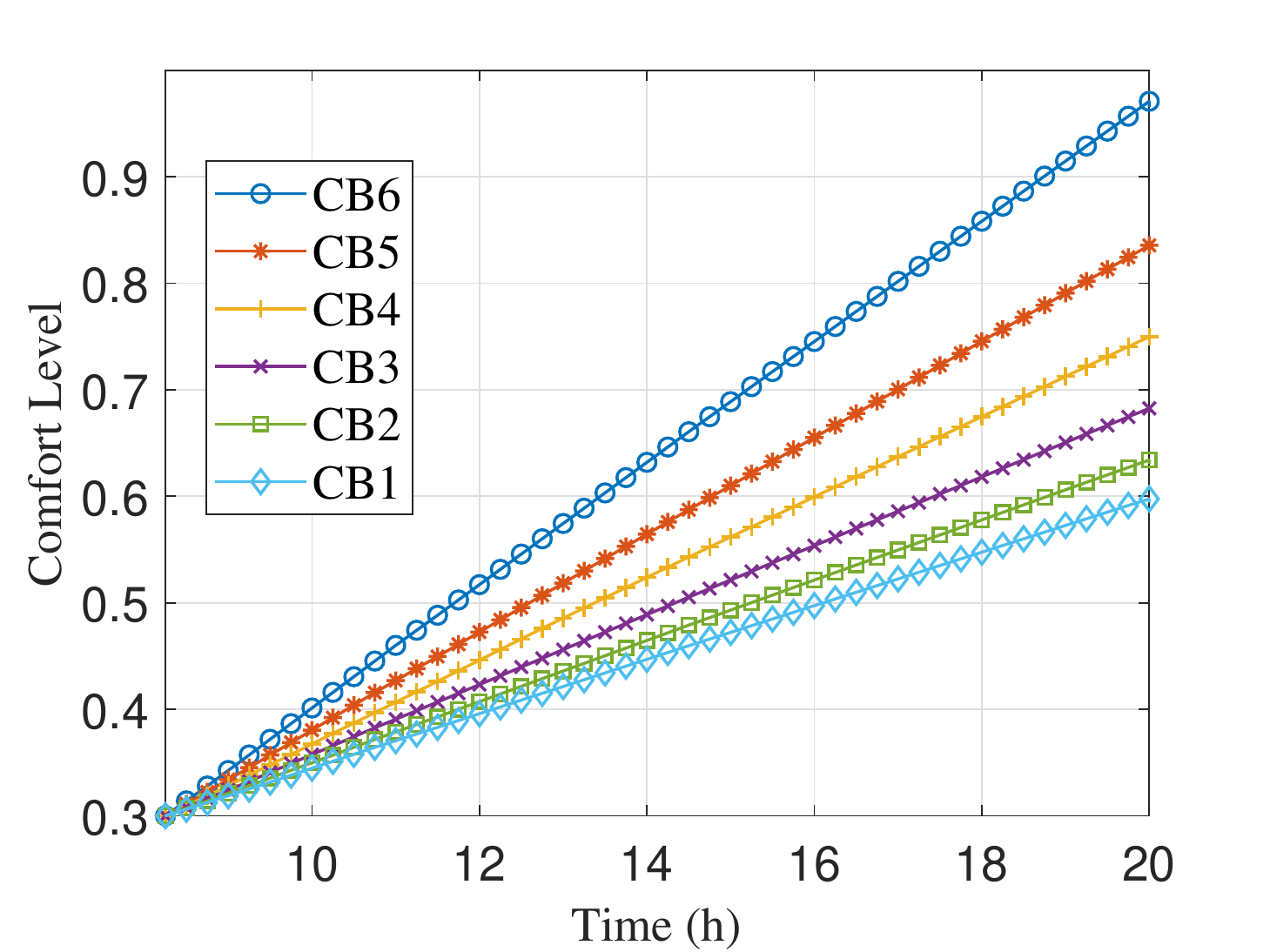}
\caption{Comfort levels related with EWHs considering our proposed comfort level model.}
\label{fig:EWH_D}
\end{figure}
\begin{figure}[!ht]
\centering
\includegraphics[width=4.8in]{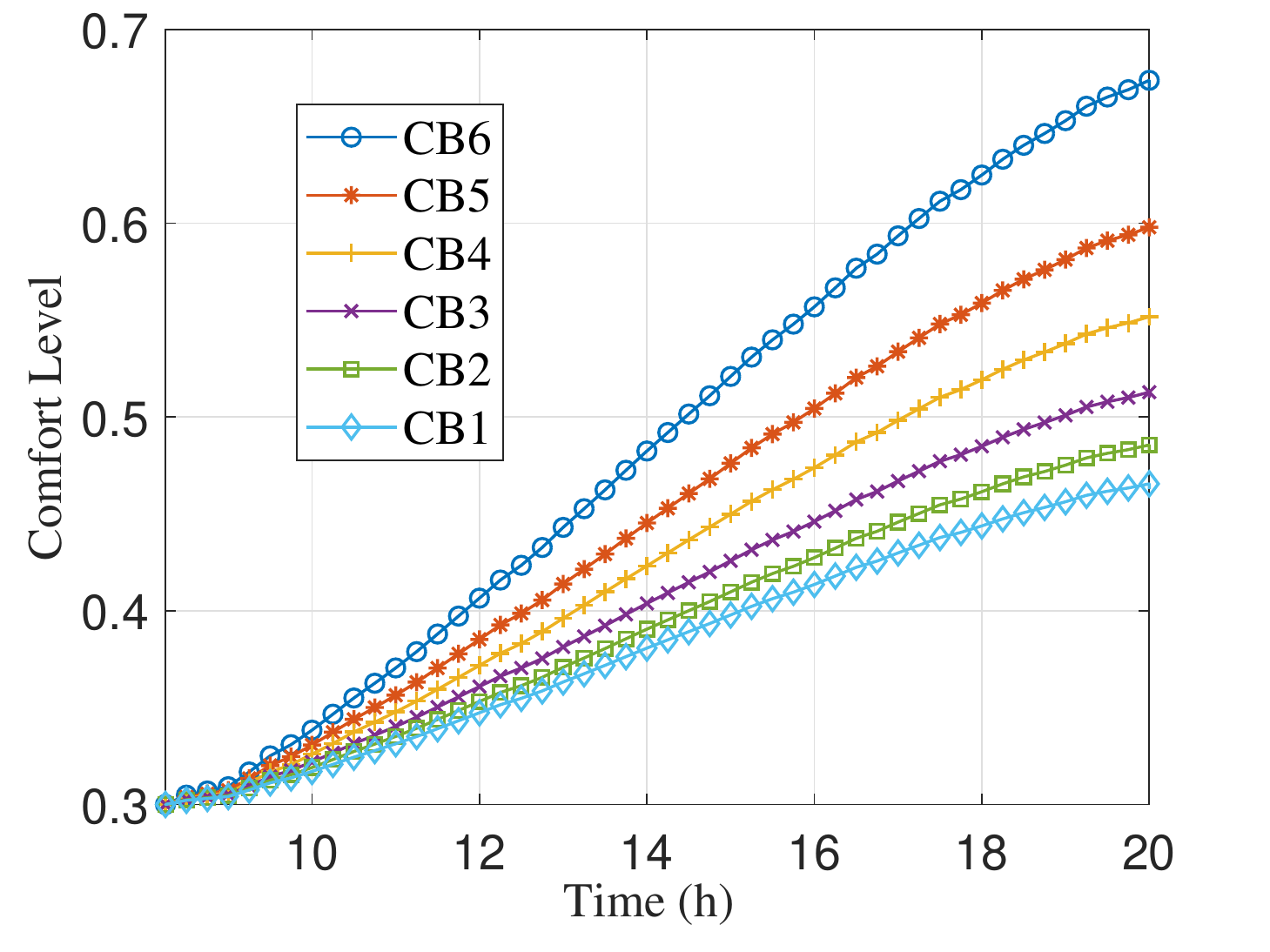}
\caption{Comfort levels related with EWHs without our proposed comfort level model.}
\label{fig:EWH_DnoJ}
\end{figure}
Moreover, the influences of comfort levels on EWH is tested. As shown in Figs.~\ref{fig:EWH_D} and~\ref{fig:EWH_DnoJ}, the comfort levels increase faster when considering the proposed comfort level models. This is because the EWHs provide $1/2$ of the total hot water demand while the rest of hot water demand is satisfied by the gas-based heat boilers. Note that the gas is far less expensive than the electricity. Besides, through checking the trade-off between electricity prices and comfort levels, comfort levels of the case with the proposed comfort level models can still maintain the comfort levels of EWHs at a high degree compared with the one without the proposed comfort level models.

\begin{table}[!ht]
\centering
\caption{Total Operating Cost and Comfort Levels of the Commercial Campus}
\label{table:5}
\begin{tabular}{|c| |c| |c| |c|}
\hline
At Noon & Without & With & $Increment\%$ \\
\hline
Total Operating Cost & 10.85 & 13.01 & $19.9\%$ \\
\hline
HVAC Comfort Level & 0.05 & 1 & $1900\%$\\
\hline
PEVs Comfort Level & 0.15-0.55 & 0.55-1 & $82\%$-$267\%$ \\
\hline
EWHs Comfort Level & 0.33-0.39 & 0.37-0.5 & $12\%$-$28\%$ \\
\hline
\end{tabular}
\end{table}
Furthermore, as shown in Table~\ref{table:5}, the pure operating cost of the one considering comfort levels is $13.01$ compared with $10.85$ of the other one without the comfort levels, which is only $19.9\%$ of increase in the operating costs but the comfort levels rise significantly as the aforementioned for one operating day. Additionally, rates of achieving the desired comfort levels are tested. Take $12:00$ pm for example, it can be observed that the comfort levels of HVAC systems increase from $5\%$ to $100\%$, which has an increment of $1900\%$; the comfort levels of PEVs increase from $15-55\%$ to $55-100\%$, which has an increment of $82-267\%$; and the comfort levels of EWHs increase from $33-39\%$ to $37-50\%$, which has an increment of $12-28\%$.

From the above results, the proposed comfort level models for HVAC systems, PEVs and EWHs work very well for a commercial campus with CBs in optimizing both operation costs and comfort levels.

\section{Conclusion}\label{sec:conclusion}
In this paper, we propose an optimal demand response framework to enable local control of demand-side appliances that are usually too small to participate directly in a wholesale electricity market. Small demand side appliances such as heating, ventilation, and air-conditioning system, electric water heaters, and plug-in electric vehicles can be utilized together by a central controller of commercial buildings to minimize the operating costs and maximize the comprehensive comfort levels. Stochastic programming has been adopted to handle various uncertainties, i.e., (i) output of renewables; (ii) day-ahead and real-time electricity prices; (iii) arrival and departure of plug-in electric vehicles; (iv) business hour demand response signals and (v) flexible energy demand. Extensive simulation results with real-world data sets indicate that our proposed model can minimize operating costs of commercial buildings while maximizing consumers' comfort levels. By sacrificing less than $20\%$ of the total operating cost, comprehensive comfort levels can increase from $0.4$ to $1$ through our proposed comfort level models.

One direction of our future work is to include vehicle to grid technology into our proposed model, i.e., the electric vehicles are treated as distributed energy storage systems. In addition, the network as well as reactive power constraints will be incorporated in our model. Finally, our proposed energy scheduling model will be verified by hardware-in-the-loop experiments.

\bibliographystyle{elsarticle-num}
\bibliography{./mybibfile}

\begin{thebibliography}{10}
\expandafter\ifx\csname url\endcsname\relax
  \def\url#1{\texttt{#1}}\fi
\expandafter\ifx\csname urlprefix\endcsname\relax\def\urlprefix{URL }\fi
\expandafter\ifx\csname href\endcsname\relax
  \def\href#1#2{#2} \def\path#1{#1}\fi

\bibitem{SaOr15}
M.~R. Sarker, M.~A. Ortega-Vazquez, D.~S. Kirschen, Optimal coordination and
  scheduling of demand response via monetary incentives, IEEE Transactions on
  Smart Grid 6~(3) (2015) 1341--1352.

\bibitem{LiAl17}
Z.~Liang, Q.~Alsafasfeh, T.~Jin, H.~Pourbabak, W.~Su, Risk-constrained optimal
  energy management for virtual power plants considering correlated demand
  response, IEEE Transactions on Smart Grid, in press.

\bibitem{BiPi15}
D.~Bian, M.~Pipattanasomporn, S.~Rahman, A human expert-based approach to
  electrical peak demand management, IEEE Transactions on Power Delivery 30~(3)
  (2015) 1119--1127.

\bibitem{KoBa16}
C.~D. Korkas, S.~Baldi, I.~Michailidis, E.~B. Kosmatopoulos, Occupancy-based
  demand response and thermal comfort optimization in microgrids with renewable
  energy sources and energy storage, Applied Energy 163 (2016) 93--104.

\bibitem{Ra16}
M.~Rahmani-andebili, Nonlinear demand response programs for residential
  customers with nonlinear behavioral models, Energy and Buildings 119 (2016)
  352--362.

\bibitem{KlKw12}
L.~Klein, J.-y. Kwak, G.~Kavulya, F.~Jazizadeh, B.~Becerik-Gerber,
  P.~Varakantham, M.~Tambe, Coordinating occupant behavior for building energy
  and comfort management using multi-agent systems, Automation in construction
  22 (2012) 525--536.

\bibitem{CB}
Energy {C}onsumption by {S}ector,
  \url{https://www.eia.gov/totalenergy/data/monthly/pdf/sec2_3.pdf} (2018).

\bibitem{Ra162}
M.~Rahmani-andebili, Modeling nonlinear incentive-based and price-based demand
  response programs and implementing on real power markets, Electric Power
  Systems Research 132 (2016) 115--124.

\bibitem{Ra15}
M.~Rahmani-andebili, Risk-cost-based generation scheduling smartly mixed with
  reliability-driven and market-driven demand response measures, International
  Transactions on Electrical Energy Systems 25~(6) (2015) 994--1007.

\bibitem{NgLe14Comfort}
D.~T. Nguyen, L.~B. Le, Joint optimization of electric vehicle and home energy
  scheduling considering user comfort preference, IEEE Transactions on Smart
  Grid 5~(1) (2014) 188--199.

\bibitem{HuRh17}
L.~Hurtado, J.~Rhodes, P.~Nguyen, I.~Kamphuis, M.~Webber, Quantifying demand
  flexibility based on structural thermal storage and comfort management of
  non-residential buildings: A comparison between hot and cold climate zones,
  Applied energy 195 (2017) 1047--1054.

\bibitem{ChXu17}
Y.~Chen, P.~Xu, Y.~Chu, W.~Li, Y.~Wu, L.~Ni, Y.~Bao, K.~Wang, Short-term
  electrical load forecasting using the support vector regression ({SVR}) model
  to calculate the demand response baseline for office buildings, Applied
  Energy 195 (2017) 659--670.

\bibitem{CuGa17}
B.~Cui, D.-c. Gao, F.~Xiao, S.~Wang, Model-based optimal design of active cool
  thermal energy storage for maximal life-cycle cost saving from demand
  management in commercial buildings, Applied Energy 201 (2017) 382--396.

\bibitem{KiNo17}
Y.~Kim, L.~K. Norford, Optimal use of thermal energy storage resources in
  commercial buildings through price-based demand response considering
  distribution network operation, Applied energy 193 (2017) 308--324.

\bibitem{HaSh18}
N.~Hajibandeh, M.~Shafie-Khah, G.~J. Os{\'o}rio, J.~Aghaei, J.~P. Catal{\~a}o,
  A heuristic multi-objective multi-criteria demand response planning in a
  system with high penetration of wind power generators, Applied Energy 212
  (2018) 721--732.

\bibitem{DoMo17}
A.~Dolatabadi, B.~Mohammadi-Ivatloo, Stochastic risk-constrained scheduling of
  smart energy hub in the presence of wind power and demand response, Applied
  Thermal Engineering 123 (2017) 40--49.

\bibitem{LiSu18}
Z.~Liang, W.~Su, Game theory based bidding strategy for prosumers in a
  distribution system with a retail electricity market, IET Smart Grid 1~(3)
  (2018) 104--111.

\bibitem{AjLu17}
A.~Ajao, J.~Luo, Z.~Liang, Q.~H. Alsafasfeh, W.~Su, Intelligent home energy
  management system for distributed renewable generators, dispatchable
  residential loads and distributed energy storage devices, in: Renewable
  Energy Congress (IREC), 2017 8th International, IEEE, 2017, pp. 1--6.

\bibitem{ZoKu12}
Y.~Zong, D.~Kullmann, A.~Thavlov, O.~Gehrke, H.~W. Bindner, Application of
  model predictive control for active load management in a distributed power
  system with high wind penetration, IEEE Transactions on Smart Grid 3~(2)
  (2012) 1055--1062.

\bibitem{ZeVa18}
I.~Zenginis, J.~Vardakas, N.~Zorba, C.~Verikoukis, Performance evaluation of a
  multi-standard fast charging station for electric vehicles, IEEE transactions
  on smart grid 9~(5) (2018) 4480--4489.

\bibitem{AhCh17}
J.~Ahn, S.~Cho, D.~H. Chung, Analysis of energy and control efficiencies of
  fuzzy logic and artificial neural network technologies in the heating energy
  supply system responding to the changes of user demands, Applied energy 190
  (2017) 222--231.

\bibitem{LiCh17}
Z.~Liang, T.~Chen, H.~Pourbabak, W.~Su, Robust distributed energy resources
  management for microgrids in a retail electricity market, in: Power Symposium
  (NAPS), 2017 North American, IEEE, 2017, pp. 1--6.

\bibitem{LiGu17}
Z.~Liang, Y.~Guo, Optimal energy management for microgrids with cogeneration
  and renewable energy sources, in: Smart Grid Communications (SmartGridComm),
  2015 IEEE International Conference on, IEEE, 2015, pp. 647--652.

\bibitem{MaJo18}
J.~M{\"a}rkle-Hu{\ss}, S.~Feuerriegel, D.~Neumann, Large-scale demand response
  and its implications for spot prices, load and policies: Insights from the
  german-austrian electricity market, Applied Energy 210 (2018) 1290--1298.

\bibitem{WaLi17}
Y.~Wang, Y.~Liu, D.~S. Kirschen, Scenario reduction with submodular
  optimization, IEEE Transactions on Power Systems 32~(3) (2017) 2479--2480.

\bibitem{Wind}
National {W}ind {T}echnology {C}enter {M}2 {T}ower,
  \url{http://www.nrel.gov/midc/nwtc_m2/} (2018).

\bibitem{PJM15}
{PJM}, \url{http://www.pjm.com/markets-and-operations.aspx} (2018).

\end{thebibliography}


\end{document}